\documentclass[11pt]{article}
\usepackage{float}
\usepackage{amsmath,amssymb,amsfonts,amsthm}
\usepackage[colorlinks,linktocpage,linkcolor=blue]{hyperref}
\usepackage{cleveref}%\usepackage[foot]{amsaddr}
\usepackage{color}
\usepackage{graphicx}
\usepackage{epstopdf}          % 如需处理EPS图像
\usepackage{float}             % 浮动体控制
\usepackage{wrapfig}           % 图文混排
\usepackage{subcaption}        % 子图支持（关键！）
\usepackage[justification=centering]{caption}
\usepackage{wrapfig}
\usepackage{multirow}
\usepackage{booktabs}
\usepackage{algorithmic}

\usepackage[T1]{fontenc}
\usepackage{authblk}
\usepackage[textwidth=6.0in,textheight=9in]{geometry}
\usepackage[linesnumbered,ruled,vlined, ruled,vlined]{algorithm2e}

\usepackage[colorlinks,linktocpage,linkcolor=blue]{hyperref}
\usepackage{fancyhdr}

\usepackage[export]{adjustbox}

\newtheoremstyle{exampstyle}
  {\topsep} % Space above
  {\topsep} % Space below
  {\itshape} % Body font
  {} % Indent amount
  {\bfseries} % Theorem head font
  {.} % Punctuation after theorem head
  {.5em} % Space after theorem head
  {} % Theorem head spec (can be left empty, meaning `normal')

\theoremstyle{exampstyle}
\numberwithin{equation}{section}
\newtheorem{theorem}{Theorem}[section]
\newtheorem{lemma}{Lemma}[section]
\newtheorem{assumption}[lemma]{Assumption}
\newtheorem{remark}[lemma]{Remark}

\newtheorem{corollary}[lemma]{Corollary}

\newtheorem{example}{Example}
\allowdisplaybreaks
\usepackage{parskip}
\let\oldref\ref
\renewcommand{\ref}[1]{(\oldref{#1})}  % stupid kludge for laTeX
\renewcommand{\eqref}[1]{(\oldref{#1})}

\newcommand{\ds}{\displaystyle}

\newbox\boxaddrone \newbox\boxaddrtwo

\newcommand{\N}{\mathcal{N}}

\newcommand{\ba}{\begin{eqnarray*}}
\newcommand{\ea}{\end{eqnarray*}}
\newcommand{\be}{\begin{equation}}
\newcommand{\ee}{\end{equation}}
\newcommand{\bea}{\begin{eqnarray}}
\newcommand{\eea}{\end{eqnarray}}

\allowdisplaybreaks
\newcommand{\p}{\partial}
\newcommand{\om}{\Omega}
\newbox\boxaddrone \newbox\boxaddrtwo

\def\N+{n\in\mathbb{N}^{+}}

\def\1d{\mathcal{D}((-\Delta)^{\gamma_1+1/2})}
\def\2d{\mathcal{D}((-\Delta)^{\gamma_2+1})}

\begin{document}
\title{Analysis and approximation for inverse potential problem of space-dependent variable-exponent sub-diffusion in heterogeneous medium} %% Article title

%\author[1]{Peiyan Guo}
\author[1]{Zhiyuan Li}
\author[2]{Chunlong Sun\thanks{Author to whom any correspondence should be addressed, sunchunlong@nuaa.edu.cn}}
\author[3]{Xiangcheng Zheng}

\affil[1]{School of Mathematics and Statistics, Ningbo University, Ningbo 315211, Zhejiang, China}
\affil[2]{School of Mathematics, Nanjing University of Aeronautics and Astronautics, Nanjing 211106, Jiangsu, China}
\affil[3]{School of Mathematics, State Key Laboratory of Cryptography and Digital Economy Security, Shandong University, Jinan 250100, China}

\maketitle

\begin{abstract}
Anomalous transport in complex heterogeneous medium is often characterized by non-Markovian memory effects and spatially varying trapping properties. From the point of view of continuous-time random walks, such heterogeneity can be described by a position-dependent heavy-tailed waiting-time density $q(t\!\!\mid\!\! x)\sim C(x)t^{-1-\alpha(x)}$, where \(0<\alpha(x)<1\) is the local anomalous exponent. This microscopic interpretation naturally leads to the sub-diffusion with a space-dependent variable exponent \(\alpha(x)\). Compared with constant-exponent models, the variable-exponent fractional operators generally do not commute with the spatial Laplacian, which makes both the forward and inverse analysis substantially more delicate. This paper investigates the inverse problem of recovering a spatially dependent radiative or potential coefficient in variable-exponent sub-diffusion from terminal observation data. We formulate the inverse problem as a nonlinear fixed-point equation and prove that the associated fixed-point operator is monotone such  that the fixed-point iteration converges decreasingly to the unique solution. We further establish a Lipschitz-type conditional stability estimate for the terminal observation map, which also implies uniqueness of the inverse problem for sufficiently large terminal time. For numerical reconstruction, we develop a fully discrete finite element scheme and derive an error estimate for the inverse problem. Numerical simulations are conducted to validate the effectiveness of the method, even in the presence of noisy data.\\

%% Keywords
\noindent Keywords: 
inverse potential problem; variable-exponent sub-diffusion; uniqueness; conditional stability; fixed-point iteration; error estimate.
\\

\end{abstract}
\section{Introduction}

Sub-diffusion models have become an important tool for describing anomalous transport in complex medium \cite{DenHou,Jinbook%JinKia,LassasLiZhang2023
}. In particular, the space-dependent variable‑exponent sub-diffusion models have been adopted in several applications. For instance, such models have been used to describe the diffusion in porous and composite medium and the dynamics of pH‑responsive microgel particles in microchannels
\cite{BrunaChapman,Pieprzyk,StiSch}. They can also describe the spatially varying effective mobility induced by ligand‑receptor binding and the hopping or crawling motion of DNA‑coated colloids
\cite{MarbachHolmes,ZhengColloid}. Inside living cells, themolecular crowding, steric barriers, and organelle interactions lead to heterogeneous transport and nonuniform lysosome distributions, which can also be modeled by variable‑exponent sub-diffusion models
\cite{BressloffNewby,FedotovLysosome,HanIntracellular,WaighKorabel}.

	Mathematically, a constant exponent in the sub-diffusion model, which corresponds to a  spatially uniform memory law, may be insufficient for the medium with varying local structures and trapping properties and thus could not properly describe the aforementioned processes. From the viewpoint of continuous-time random walks, such spatial heterogeneity can be described by a position-dependent heavy-tailed waiting-time density $q(t\!\mid \!x)\sim C(x)t^{-1-\alpha(x)}$, where \(\alpha(x)\) is the local anomalous exponent \cite{Angstmann2026,CheGorSok,FedFal,FedHan,Straka}. This microscopic interpretation naturally leads to the space-dependent variable-exponent sub-diffusion.
	
	Let $\Omega \subset \mathbb{R}^d$ with $1 \le d \le 3$ be a simply-connected bounded domain with the piecewise smooth boundary $\p \om$ with convex corners, $0<T<\infty$ be the terminal time, $x := (x_1,\cdots,x_d)$ denotes the spatial variables, respectively. The variable-exponent sub-diffusion model is given by 
	\begin{equation}\label{Prob-IBVP}
		\begin{cases}
			\begin{aligned}
				\partial_t^{\alpha(x)} u - \Delta u + q(x)u &= f(x), &  &(x,t) \in \Omega_T := \Omega \times (0,T], \\
				u &= u_{0}(x),  &  &(x,t) \in \Omega\times\{0\}, \\
				u&= u_b(x), &  &(x,t) \in \partial\Omega \times (0,T],
			\end{aligned}
		\end{cases}
	\end{equation}
	where $q$, $f$, $u_0$ and $u_b$ denote the potential, the source term, the initial value and the boundary value, respectively. Throughout the paper we assume $0<\alpha_*:=\operatorname*{ess\,inf}_{x\in\Omega}\alpha(x)
	\leq \alpha(x)\leq
	\operatorname*{ess\,sup}_{x\in\Omega}\alpha(x)=:\alpha^*<1.$
	The Caputo derivative of space-dependent exponent is defined by \cite{KiaSocYam}
	\begin{equation*}%\label{m3as}
		\partial_t^{\alpha(x)}u(x,t)
		:=\frac{1}{\Gamma(1-\alpha(x))}
		\int_0^t (t-s)^{-\alpha(x)}\,\partial_s u(x,s)\,ds.
	\end{equation*}

Direct and inverse problems for sub-diffusion have been extensively studied in the constant-exponent case \cite{AmmariKang2006,DElia2020,DuYangZhou2020,DuZhou2023,JinKia,LassasLiZhang2023,LinLiu2023,zhouzhi2022}. In contrast, the space-dependent variable-exponent model brings additional difficulties, e.g. the corresponding fractional derivative generally does not commute with the spatial operator, such that only a few investigations are available. For instance, an $\alpha$-robust semi-discrete finite element method has been studied for a space-dependent variable-exponent time-fractional Fokker--Planck initial-boundary value problem in \cite{LeStynes}. The analysis of a fully-discrete numerical scheme can be found in \cite{Nie2025}. Recently, significant progresses have been achieved for inverse problems of such models, cf. the identification of the space-dependent variable exponent from boundary data \cite{Hon,Hong2026,Jin2026}.  Nevertheless, the recovery of a potential coefficient from terminal observation for variable-exponent sub-diffusion remains less explored.

Our objective is to use the final time data $u(x,T)$ to recover the unknown source term $q(x)$. For $q(x)$, we define $\mathcal{Q}:=\{\psi\in C(\overline\Omega): 0\leq \psi \leq M_1\}$ as the admissible set, and denote the forward operator $\mathbb{G}: \mathcal{Q}\to H^2(\Omega)$ by $\mathbb{G}\psi:=u(x,T;\psi),\;x\in\Omega.$
	Then the inverse problem of interest could be stated as follows: 
	\begin{equation}\label{IP-continuous}
		\text{with data}\ \phi(x):=u(x,T),\ \text{find}\ q\in {\mathcal Q}\ \text{such that}\ \mathbb{G} q=\phi.
	\end{equation}
	We also consider the numerical reconstruction from a noisy datum $\phi^\delta(x)=\phi(x)+\xi(x)$ in $\Omega$, and $\xi$ denotes the measurement noise. The accuracy of the observational data $\phi^\delta$ is measured by the noise level $\|\phi^\delta - \phi\|_{C({\bar\Omega})}= \delta$. Such an inverse problem may arise in many practical applications, where $q$ represents the radiativity coefficient in heat conduction \cite{yamamoto2001simultaneous} and the perfusion coefficient in Pennes' bio-heat equation in human physiology \cite{pennes1948analysis}.
	
	In this work, the uniqueness and stability of the nonlinear inverse problem \eqref{IP-continuous} will be investigated. Denoting the exact solution by $q^*$, we construct the fixed-point operator
	$K q=\frac{-\partial^{\alpha(x)}_t u(x,T;q)+\Delta \phi(x)+f(x)}{\phi(x)} \text{ for }  q \in \mathcal{Q},$
	whose fixed points generate the desired data $\phi(x)$; here the notation $u(x,T;q)$ reflects the dependence of $u$ on $q$. We prove that $K$ is monotone increasing (Lemma \ref{lemma_monotone}), so that the fixed-point iteration converges decreasingly to the unique solution (Theorem \ref{theorem_uniqueness}). In addition, we establish a Lipschitz-type stability estimate for the terminal observation map (Theorem \ref{thm-stability}), which also implies the uniqueness of the inverse problem for sufficiently large terminal time (Corollary \ref{cor-uniqueness}). The fixed-point iteration is then used as the numerical reconstruction scheme in Section \ref{Algorithm}. 
Furthermore, we develop a fully discrete finite element scheme and derive an error estimate for the inverse problem. We show that there exists \(T_*>0\) such that for every \(T\ge T_*\),
\begin{equation*}
\|q^*-q^\delta\|_{L^2(\Omega)}\le
		C\left(\frac{\delta}{h^2}+h+\tau\right),
\end{equation*}
where $q^*$ and $q^\delta$ are the exact solution and its reconstruction, respectively. Here $h$ and $\tau$ denote the space mesh size and time step size, respectively. This a priori error estimate provides guidelines to choose discretization parameters $h$ and $\tau$ according to the noise level $\delta$. For example, if the discretization parameters are chosen as $h\sim \tau\sim \delta^{1/3}$, then the reconstruction error satisfies $\|q^*-q^\delta\|_{L^2(\Omega)}=\mathcal O(\delta^{1/3}).$ Numerical simulations are conducted to validate the effectiveness of the method, demonstrating its ability to accurately reconstruct even in the presence of noisy data.

	The remainder of the paper is organized as follows. Section \ref{section2-directProb} collects the forward results used in the inverse analysis. Section \ref{section3-IP} studies the fixed-point formulation, uniqueness, and conditional stability. Error estimate for the inverse problem is developed in Section \ref{Algorithm}, and numerical reconstruction is described in Section \ref{sec5}. The conclusion
	is given in the last section.
	\section{Well-posedness of forward problem and positivity results}\label{section2-directProb}
	For subsequent analysis, we first recall some standard results for parabolic equations. Let $L^p(\Omega)$ with $1 \le p \le \infty$ be the Banach space of $p$th power Lebesgue integrable functions on $\om$. Denote the inner product of $L^2(\Omega)$ by $(\cdot,\cdot)$. For a positive integer $m$,
	let $ W^{m, p}(\Omega)$ be the Sobolev space of $L^p$ functions with weakly derivatives of $m$ in $L^p(\Omega)$. Let $H^m(\Omega) := W^{m,2}(\Omega)$ and $H^m_0(\Omega)$ be its subspace with the zero boundary condition up to the order $m-1$. For $\frac{1}{2}<s\leq 1$, $H^s_0(\Omega)$ is defined by interpolating $L^2(\Omega)$ and $H^1_0(\Omega)$ \cite{AdaFou}. Let $I$ be the identity operator and $\{\lambda_n,\varphi_n\}_{n=1}^\infty$ be the Dirichlet eigensystem of the operator $-\Delta+q I$, where $0<\lambda_1\leq \lambda_2\leq \cdots$ and the corresponding eigenfunctions $\{\varphi_n\}_{n=1}^\infty$ form an orthonormal basis of $L^2(\Omega)$.
	
	For a Banach space $X$, we say that $f\in L^p(0,T;X)$ provided that
	\begin{eqnarray*}%\label{Lpde}
		\|{f}\|_{L^p(0,T;X)}:=
		\begin{cases}
			\ds\left(\int_0^T \|{f(\cdot,t)}\|_X^p {\rm d}t\right)^{1/p} <\infty,&  1\leq p <\infty,\\[0.1in]
			\ds{\rm ess\; sup}_{0<t<T} \|{f(\cdot,t)}\|_X <\infty, &  p=\infty.
		\end{cases}
	\end{eqnarray*}
	Let $W^{m, p}(0,T; X)$ be the space of functions in $W^{m, p}(0,T)$ with respect to $\|\cdot\|_X$.  Similarly, for $0\leq t_0 <T$, let $C([t_0,T];X)$ be the space of functions in $C([t_0,T])$ with respect to $\|\cdot\|_X$. Moreover, define $C((0,T];X):=\cap_{0<t_0<T} C([t_0,T];X)$. All spaces are equipped with standard norms \cite{AdaFou,Eva}. 
	\subsection{Well-posedness of forward problem}
	We recall the admissible set of unknown $q$ as $\mathcal{Q}:=\{\psi\in C(\overline\Omega): 0\leq \psi \leq M_1\}$.
	First, we give the well-posedness of the forward problem \eqref{Prob-IBVP}.

	\begin{lemma}\label{wellposedness-ivp}
		Assume that $q\in\mathcal Q$, that the Dirichlet realization of $-\Delta$ has the usual $H^2$-regularity, and that $u_b$ is in the trace space $H^{\frac32}(\partial\Omega)$. If $f\in L^2(\Omega)$ and $u_0\in L^2(\Omega)$, then the problem \eqref{Prob-IBVP} has a unique weak solution $u\in C((0,T];L^2(\Omega))$. Moreover, for every $\tau\in(0,T)$, we have $u\in C([\tau,T];H^2(\Omega)), 
		\partial_t^{\alpha(x)}u\in C([\tau,T];L^2(\Omega))$. If $u_0\in H^2(\Omega)$ and $u_0|_{\partial\Omega}=u_b$, the initial and boundary data are compatible in the usual trace sense.
	\end{lemma}
	
	\begin{proof}
		Let $A_q$ be the Dirichlet realization of $-\Delta+q$ in $L^2(\Omega)$, with domain
		$
		D(A_q)=H^2(\Omega)\cap H_0^1(\Omega).$
		The graph norms of $A_q$, $q\in\mathcal Q$, are uniformly equivalent to the $H^2(\Omega)$-norm. Let $U_q$ denote the stationary solution $-\Delta U_q+qU_q=f~\text{in }\Omega,
		U_q=u_b~\text{on }\partial\Omega.$
		Choose \(b\in H^2(\Omega)\) such that \(b|_{\partial\Omega}=u_b\), whose existence follows from the trace theorem since \(u_b\in H^{3/2}(\partial\Omega)\). Writing \(U_q = b + \widetilde{U}_q\) reduces this problem to $A_q\widetilde U_q=f+\Delta b-qb, \widetilde U_q\in D(A_q).$
		By the Lax--Milgram theorem and elliptic regularity results, there exists a unique $U_q\in H^2(\Omega)$.
		Set $v=u-U_q$ and $h_q=u_0-U_q$. Then $v$ satisfies the homogeneous problem $\partial_t^{\alpha(x)}v+A_qv=0, v|_{\partial\Omega}=0, v(0)=h_q.$
		
		For $z\in\mathbb C\setminus(-\infty,0]$, let $M_z$ be the multiplication operator by $z^{\alpha(x)}$, with the principal branch, and put $R_q(z)=(A_q+M_z)^{-1}.$ By \cite[Proposition~2.1 and Theorem~1.1]{KiaSocYam}, $R_q$ is holomorphic on the cut plane and maps $L^2(\Omega)$ into $D(A_q)$. For $\varepsilon\in(0,1)$ and $\theta\in(\pi/2,\pi)$, define the contour $\gamma(\varepsilon,\theta)=\gamma_-(\varepsilon,\theta)\cup\gamma_0(\varepsilon,\theta)\cup\gamma_+(\varepsilon,\theta),$ where $\gamma_\pm(\varepsilon,\theta)=\{re^{\pm i\theta}:r\geq\varepsilon\},
		\gamma_0(\varepsilon,\theta)=\{\varepsilon e^{i\beta}:-\theta\leq\beta\leq\theta\}.$
		The variable-exponent resolvent theorem in \cite[Theorem~1.1]{KiaSocYam} shows that the unique mild solution of the homogeneous problem is given by
		\begin{equation}\label{eq:forward-contour}
			v(t)=\frac{1}{2\pi i}\int_{\gamma(\varepsilon,\theta)}
			e^{zt}R_q(z)z^{-1}M_z h_q\,dz.
		\end{equation}
		The integral is independent of $\varepsilon\in(0,1)$. Because the contour avoids the negative real axis and $R_q$ is holomorphic, the mapping $t\mapsto v(t)$ is analytic from $(0,\infty)$ into $L^2(\Omega)$. Hence $v\in C^\infty((0,\infty);L^2(\Omega))$, and in particular $v'\in C((0,\infty);L^2(\Omega))$.
		
		We now verify that $v$ is a strong solution. From the definition of the Caputo derivative and the analyticity of $v$, 
		$
		\partial_t^{\alpha(x)}v(t)
		=\frac{1}{\Gamma(1-\alpha(x))}
		\int_0^t (t-s)^{-\alpha(x)}v'(s)\,ds.
		$
		For $t\geq\tau>0$, the integral converges in $L^2(\Omega)$ uniformly in $t$ on compact subsets of $(0,\infty)$, since $v'\in C([\tau/2,T];L^2(\Omega))$. Thus $\partial_t^{\alpha(x)}v\in C([\tau,T];L^2(\Omega))$. Moreover, the contour representation \eqref{eq:forward-contour} gives $v(t)=\mathcal O(t^{\alpha_*})$ in $L^2(\Omega)$ as $t\to0^+$, hence $\partial_t^{\alpha(x)}v$ is locally integrable and of exponential type on $(0,\infty)$; its Laplace transform is therefore well defined.
		
		To identify $\partial_t^{\alpha(x)}v$, we use the resolvent equation. Since $v$ is the mild solution, its Laplace transform satisfies $\widehat v(z)=R_q(z)z^{-1}M_z h_q$ for $z$ in a suitable sector. Taking the Laplace transform of $\partial_t^{\alpha(x)}v$ and using the Caputo-derivative rule
		$\widehat{\partial_t^{\alpha(x)}v}(z)=M_z\widehat v(z)-z^{-1}M_z h_q$ (applied pointwise in $x$; valid for sufficiently smooth $v$), we get
		\[
		\widehat{\partial_t^{\alpha(x)}v}(z)
		=M_z\widehat v(z)-z^{-1}M_z h_q
		=M_z R_q(z)z^{-1}M_z h_q - z^{-1}M_z h_q.
		\]
		On the other hand, from $R_q(z)=(A_q+M_z)^{-1}$ we have $A_qR_q(z)=I-M_zR_q(z).$ Therefore,
		\[
		\widehat{\partial_t^{\alpha(x)}v}(z)
		=(I-A_qR_q(z))z^{-1}M_z h_q - z^{-1}M_z h_q
		=-A_qR_q(z)z^{-1}M_z h_q
		=-A_q\widehat v(z).
		\]
		By uniqueness of the Laplace transform (or by direct inversion along the contour), we conclude that $\partial_t^{\alpha(x)}v(t)=-A_qv(t) \text{in } L^2(\Omega)\text{ for every } t>0.$ Hence $\partial_t^{\alpha(x)}v=-A_qv\in C([\tau,T];L^2(\Omega))$.
		
		Finally, since $u=U_q+v$, we obtain $u\in C([\tau,T];H^2(\Omega)),$ and $\partial_t^{\alpha(x)}u=\partial_t^{\alpha(x)}v\in C([\tau,T];L^2(\Omega)),$ proving the asserted regularity. The uniqueness of the weak solution follows from the uniqueness of the mild solution, which in turn is a consequence of the uniqueness of the Laplace transform.
	\end{proof}
	
	\subsection{Maximum principle and positive results}
	
	Next, we introduce the maximum principle for the diffusion model \eqref{Prob-IBVP}. 
	
	\begin{lemma}\label{lemma_maximum}
		Let $u$ be a classical solution of \eqref{Prob-IBVP} such that $u(x,\cdot)\in W^{1,1}(0,T)$ for every $x\in\overline\Omega$. If $q$, $f$, $u_0$ and $u_b$ are nonnegative, then $u(x,t)\geq0,\qquad (x,t)\in\overline\Omega\times[0,T].$ The same conclusion holds almost everywhere for weak solutions obtained as $L^2$-limits of nonnegative smooth approximations.
	\end{lemma}
	\begin{proof}
		Suppose, to the contrary, that $u$ takes a negative value. By continuity, its minimum on $\overline\Omega\times[0,T]$ is attained at some $(x_0,t_0)$. The nonnegative initial and boundary data imply that $x_0\in\Omega$, $t_0>0$, and $u(x_0,t_0)<u(x_0,0).$ Put $\beta=\alpha(x_0)$. The extremum identity for the Caputo derivative gives \cite[Theorem~1]{yury2009}
		\begin{align*}
			\partial_t^\beta u(x_0,t_0)
			=\frac{u(x_0,t_0)-u(x_0,0)}{\Gamma(1-\beta)t_0^\beta} +\frac{\beta}{\Gamma(1-\beta)}
			\int_0^{t_0}\frac{u(x_0,t_0)-u(x_0,s)}{(t_0-s)^{1+\beta}}\,ds<0.
		\end{align*}
		At the spatial minimum, $\Delta u(x_0,t_0)\geq0$, while $q(x_0)u(x_0,t_0)\leq0$, $\partial_t^{\alpha(x_0)}u(x_0,t_0)-\Delta u(x_0,t_0)
		+q(x_0)u(x_0,t_0)<0,$ which contradicts $f(x_0)\geq0$. Therefore $u$ is nonnegative.
		
		It remains to justify the passage to weak solutions. Let $u$ be the weak solution provided by Lemma \ref{wellposedness-ivp} for nonnegative data $q,f,u_0,u_b$ as in the statement. Approximate the data by smooth nonnegative functions \(q_\epsilon\), \(f_\epsilon\), \(u_{0,\epsilon}\), \(u_{b,\epsilon}\) with \(u_{0,\epsilon}=u_{b,\epsilon}\) on \(\partial\Omega\) and \(u_{b,\epsilon}\to u_b\) in \(H^{3/2}(\partial\Omega)\) (for instance by mollification, retaining the trace compatibility).
		For smooth data the corresponding solution $u_\varepsilon$ is classical: it is analytic in $t$ for $t>0$ by the contour representation of Lemma \ref{wellposedness-ivp}, and elliptic regularity in $x$ supplies the spatial smoothness needed to run the pointwise argument above, so $u_\varepsilon\ge0$. Moreover the solution map is continuous in the data: 
		$q\mapsto R_q(z)=(A_q+M_z)^{-1}$ is continuous in the operator norm by the resolvent identity $R_{q_1}(z)-R_{q_2}(z)=-R_{q_1}(z)(q_1-q_2)R_{q_2}(z),$ and \(h_q = u_0 - U_q\) depends continuously on \((f, u_0, u_b)\), whence \eqref{eq:forward-contour} gives $u_\varepsilon\to u \text{ in } C((0,T];L^2(\Omega)),\quad \text{ as } \varepsilon\to 0.$ Passing to a subsequence yields $u_\varepsilon(x,t)\to u(x,t)$ for a.e.\ $(x,t)$, and therefore $u\ge0$ a.e.\ in $\Omega\times(0,T]$.
	\end{proof}
	
	Throughout the paper, we need the following assumptions to be valid.
	\begin{assumption}\label{assumption}
		The exact potential $q^*$, the initial data $u_{0}$, the boundary data $u_b$ and the source term $f$ are continuous and satisfy the following conditions: $(a)~q^*\in \mathcal{Q}:=\{q\in C(\bar\Omega): 0\leq q \leq M_1\}$, where $M_1>0$ is the upper bound of the admissible set $\mathcal Q$; $(b)~u_{0}\geq M_{2}>0 \mathrm{~in~} \Omega$, $u_0(x)=u_b(x)$ for all $x\in \partial\Omega$; 
		$(c)~u_b \geq M_2>0$ on $\partial\Omega$; $(d)~f\geq 0$ and $f\geq -\Delta u_0+M_{1}u_0 \ \mathrm{in} \ \Omega$.
	\end{assumption}

	With Assumption \ref{assumption} and Lemma \ref{lemma_maximum}, we could deduce the positivity results of the solutions of equations \eqref{Prob-IBVP}.

	\begin{lemma}\label{lemma_positivity}
		With Assumption \ref{assumption}, the \(u\) defined in model \eqref{Prob-IBVP} satisfies the following results for every \(q\in\mathcal Q\): $(a)~u(x,t)\ge u_0(x)\ge M_2$ for a.e. \((x,t)\in \Omega\times(0,T]\); $(b)~\partial_t^{\alpha(x)}u(x,t) \ge 0$ for a.e. \((x,t)\in \Omega\times(0,T]\), and
		\(f(x)+\Delta u(x,t) \ge q(x)M_2\) for a.e. \(x\in \Omega\) and every \(t>0\).
	\end{lemma}
	\begin{proof}
		For statement (a), set \(w=u-u_0\). Since \(u_0\) is time-independent and \(u_0=u_b\) on \(\partial\Omega\) by Assumption \ref{assumption}(b), \(w\) satisfies
		$$\partial_t^{\alpha(x)} w - \Delta w + q(x)w = f(x)+\Delta u_0(x)-q(x)u_0(x), (x,t) \in \Omega_T,$$
		equipped with homogeneous initial and boundary conditions.
		%\begin{equation*}
		%\begin{cases}
		%\begin{aligned}
		%	\partial_t^{\alpha(x)} w - \Delta w + q(x)w &= f(x)+\Delta u_0(x)-q(x)u_0(x), &  &(x,t) \in \Omega_T, \\
		%w &= 0,  &  &(x,t) \in \partial\Omega \times (0,T], \\
		%	w&= 0, &  &(x,t) \in \Omega\times\{0\}.
		%\end{aligned}
		%\end{cases}
		%\end{equation*}
		By Assumption \ref{assumption}(d) and \(q\le M_1\) we have
		\(f(x)+\Delta u_0(x)-q(x) u_0(x) \ge f(x)+\Delta u_0(x)-M_1u_0(x) \ge 0\) a.e. in $\Omega$.
		Lemma \ref{lemma_maximum} therefore gives \(w\ge0\) a.e. in \(\Omega\times(0,T]\), i.e.\ \(u(x,t)\ge u_0(x)\ge M_2\) for a.e.\ \((x,t)\in\Omega\times(0,T]\), which proves (a).
		
		For statement (b), let \(A_q\) be the Dirichlet realization of \(-\Delta+q\) in \(L^2(\Omega)\) introduced in the proof of Lemma \ref{wellposedness-ivp}, write \(u=u_0+z\), and put \(g=f+\Delta u_0-qu_0\), which is nonnegative a.e.\ by Assumption \ref{assumption}(d) and \(q\le M_1\). Then \(z\) solves $\partial_t^{\alpha(x)}z+A_qz=g, z|_{\partial\Omega}=0, z(0)=0.$
		Let \(M_s\) denote the multiplication operator by \(s^{\alpha(x)}\) and put \(G(s)=(A_q+M_s)^{-1}\) for \(s>0\). By the resolvent representation underlying Lemma \ref{wellposedness-ivp}, the mild solution satisfies $z(t)=\mathcal L^{-1}\big[G(s)g/s\big](t),$ where \(\mathcal L^{-1}\) denotes the inverse Laplace transform and \(g\) is regarded as a time-independent function.
		
		We claim that for all nonnegative \(h,\widetilde h\in L^2(\Omega)\), the scalar function $s\;\longmapsto\;\langle G(s)h,\widetilde h\rangle$ is completely monotone on \((0,\infty)\). To see this, differentiate the resolvent equation: by induction on \(k\),
		\[
		G^{(k)}(s)=\sum_{j=1}^{k}(-1)^j\sum_{\ell_1+\cdots+\ell_j=k,\ \ell_i\ge1}
		\frac{k!}{\ell_1!\cdots\ell_j!}\,
		G(s)M_s^{(\ell_1)}G(s)\cdots M_s^{(\ell_j)}G(s),
		\]
		where \(M_s^{(\ell)}\) is the multiplication by
		\(\alpha(x)(\alpha(x)-1)\cdots(\alpha(x)-\ell+1)\,s^{\alpha(x)-\ell}\).
		Since \(0<\alpha(x)<1\), the factor \((-1)^{\ell-1}\alpha(\alpha-1)\cdots(\alpha-\ell+1)\) is nonnegative for every \(\ell\ge1\), so that $\widetilde M_\ell:=(-1)^{\ell-1}M_s^{(\ell)}\ge0$ as a multiplication operator. Using \(k-j=\sum_{i=1}^j(\ell_i-1)\), we may rewrite
		$
		(-1)^kG^{(k)}(s)=\sum_{j,\ell}\frac{k!}{\ell_1!\cdots\ell_j!}\,
		G(s)\widetilde M_{\ell_1}G(s)\cdots\widetilde M_{\ell_j}G(s).
		$
		Moreover, \(G(s)\) is positivity preserving: for \(f\ge0\), the function \(u=G(s)f\) solves
		\(-\Delta u+(q+s^{\alpha(x)})u=f\) in \(\Omega\) with \(u=0\) on \(\partial\Omega\), and the elliptic maximum principle gives \(u\ge0\). Consequently every summand maps nonnegative functions to nonnegative functions, and for \(h,\widetilde h\ge0\),
		\[
		(-1)^k\frac{d^k}{ds^k}\langle G(s)h,\widetilde h\rangle
		=\big\langle (-1)^kG^{(k)}(s)h,\widetilde h\big\rangle\ge0,
		\qquad k=0,1,2,\dots,
		\]
		which proves the complete monotonicity. Bernstein's theorem for completely monotone functions (see, e.g., \cite{SSV}) now provides, for every nonnegative \(h,\widetilde h\in L^2(\Omega)\), a nonnegative finite measure \(\nu=\nu_{h,\widetilde h}\) on \([0,\infty)\) such that $\langle G(s)h,\widetilde h\rangle=\int_0^\infty e^{-st}\,d\nu(t).$
		Polarizing this identity over the positive cone of \(L^2(\Omega)\) extends it to all \(h,\widetilde h\in L^2(\Omega)\), with \(h\mapsto\nu_{h,\widetilde h}\) linear in \(h\); this bilinear-form (equivalently, operator-valued) form of Bernstein's theorem is what is used below.
		
		Since \(1/s\) is the Laplace transform of the Heaviside function \(H\), taking \(h=g\) and an arbitrary nonnegative \(\widetilde h\in L^2(\Omega)\) yields, by the convolution theorem for the Laplace transform,
		\[
		\langle z(t),\widetilde h\rangle
		\!=\!\mathcal L^{-1}\big[\langle G(s)g,\widetilde h\rangle/s\big](t)
		\!=\!(\nu_{g,\widetilde h}*H)(t)
		\!=\!\int_{[0,t]}H(t-r)\,d\nu_{g,\widetilde h}(r)
		\!=\!\nu_{g,\widetilde h}\big([0,t]\big)\!\ge\!0,
		\]
		and \(t\mapsto\langle z(t),\widetilde h\rangle=\nu_{g,\widetilde h}([0,t])\) is nondecreasing. Hence \(z(t)\ge0\) for every \(t\in(0,T]\). 
		Since $z \in C^\infty((0,\infty); L^2(\Omega))$ by Lemma 2.1, the map $t \mapsto z(t)$ is continuous in $L^2(\Omega)$. Combining this continuity with the monotonicity of $t \mapsto \langle z(t), \widetilde{h} \rangle$ for all nonnegative $\widetilde{h} \in L^2(\Omega)$, we claim that there exists a representative of $z$ such that $z(x, t_2) \ge z(x, t_1)$ for a.e. $x \in \Omega$ and all $0 \le t_1 < t_2 \le T$. Indeed, by the separability of $L^2(\Omega)$, first establish the monotonicity on a countable dense set of rational times in $[0,T]$ (where the null set can be chosen uniformly thanks to the $L^2$-continuity of $t \mapsto z(t)$), and then extend by continuity to all times. Since $z$ is smooth in $t$ for $t > 0$ by Lemma 2.1, this representative is smooth (hence absolutely continuous) in $t$ on $(0,T]$ for a.e. $x \in \Omega$, with $\partial_t z(x,t) \ge 0$ a.e.
		\ Since \(\partial_t^{\alpha(x)}u_0=0\) and \(z(x,0)=0\), the definition of the Caputo derivative gives, for a.e.\ \((x,t)\in\Omega\times(0,T]\),
		\[
		\partial_t^{\alpha(x)}u(x,t)=\partial_t^{\alpha(x)}z(x,t)
		=\frac1{\Gamma(1-\alpha(x))}\int_0^t (t-s)^{-\alpha(x)}\partial_s z(x,s)\,ds\ge0.
		\]
		This together with equation \eqref{Prob-IBVP}, part (a) and \(q\ge0\) yields
		\[
		f(x)+\Delta u(x,t)= \partial_t^{\alpha(x)} u(x,t) + q(x) u(x,t) \ge q(x) u(x,t) \ge q(x) M_2 \;\;  {\rm a.e.} \;\; \Omega.
		\]
		The proof is complete.
	\end{proof}

	\subsection{Time uniform-bounds of solution}
	
	\begin{lemma}\label{lem:frac_deriv_estimate}
		Let the assumptions of Lemma \ref{wellposedness-ivp} hold and let $U_q$ be the stationary solution introduced there. Then, for every $\tau>0$ and every $T>\tau$, there is a constant $C_\tau>0$, independent of $q\in\mathcal Q$ and $T$, such that
		\begin{equation}\label{eq:uniform-terminal-bounds}
			\sup_{\tau\leq t\leq T}\left(
			\|u(\cdot,t;q)\|_{H^2(\Omega)}
			+\|\partial_t^{\alpha(x)}u(\cdot,t;q)\|_{L^2(\Omega)}
			\right)\leq C_\tau.
		\end{equation}
		Moreover, there is a constant $C>0$, independent of $q$ and $t$, for which
		\begin{equation}\label{eq:stationary-relaxation}
			\|u(\cdot,t;q)-U_q\|_{H^2(\Omega)}
			+\|\partial_t^{\alpha(x)}u(\cdot,t;q)\|_{L^2(\Omega)}
			\leq Ct^{-\alpha_*},\qquad t\geq1.
		\end{equation}
	\end{lemma}
	\begin{proof}
		Write $v_q=u(\cdot,\cdot;q)-U_q$ and $h_q=u_0-U_q$. Uniform elliptic regularity gives $\sup_{q\in\mathcal Q}\bigl(\|U_q\|_{H^2(\Omega)}+\|h_q\|_{L^2(\Omega)}\bigr)<\infty.$ With the notation of \eqref{eq:forward-contour},
		\begin{equation}\label{eq:vq-contour}
			v_q(t)=\frac{1}{2\pi i}\int_{\gamma(\varepsilon,\theta)}
			e^{zt}R_q(z)z^{-1}M_z h_q\,dz.
		\end{equation}
		Proposition~2.1 of \cite{KiaSocYam} states that $R_q(z)=(A_q+M_z)^{-1}$ maps $L^2(\Omega)$ into $D(A_q)$ and that, on the contour $\gamma(\varepsilon,\theta)$ (on which the constant $C(r,\beta)$ in \cite[(2.5)]{KiaSocYam} depends only on $\alpha_*$, $\alpha^*$ and $\theta$, hence is independent of $q$), $\|R_q(z)\|_{L^2\to L^2}\le C\max\big(|z|^{-\alpha_*},|z|^{-\alpha^*}\big).$
		Combining this with $A_qR_q(z)=I-M_zR_q(z)$ and $\|M_z\|_{L^2\to L^2}\le\max(|z|^{\alpha_*},|z|^{\alpha^*})$ yields constants $r_0\in(0,1)$ and $C>0$, uniform in $q$, such that
		\begin{align}
			\|R_q(z)\|_{L^2\to D(A_q)}+
			\|A_qR_q(z)\|_{L^2\to L^2}&\leq C,
			&&0<|z|\leq1,\label{eq:small-resolvent}\\
			\|R_q(z)\|_{L^2\to D(A_q)}+
			\|A_qR_q(z)\|_{L^2\to L^2}&\leq C(1+|z|^{\alpha^*-\alpha_*}),
			&&|z|\geq1.\label{eq:large-resolvent}
		\end{align}
		The second estimate follows from $\|R_q(z)\|_{L^2\to L^2}\!\le\! C|z|^{-\alpha_*}, |z|\!\ge\!1$ and $\|M_z\|_{L^2\to L^2}\!\le\!|z|^{\alpha^*}$, which via $A_qR_q(z)=I-M_zR_q(z)$ give $\|A_qR_q(z)\|_{L^2\to L^2}\le1+C|z|^{\alpha^*-\alpha_*}$. For small $|z|$, the first estimate follows first from the Neumann representation $(A_q+M_z)^{-1}=(I+A_q^{-1}M_z)^{-1}A_q^{-1}$ when $|z|\leq r_0$ (here $q\ge0$ gives $\|A_q^{-1}\|\le\|(-\Delta)^{-1}\|$, uniformly in $q$), and then from the quoted resolvent result on the compact annulus $r_0\leq|z|\leq1$ (on which the $L^2\to L^2$ estimate above is uniform in $q$, and $\|M_z\|_{L^2\to L^2}\le1$ together with $A_qR_q(z)=I-M_zR_q(z)$ bound $A_qR_q(z)$). Notice that $\|z^{-1}M_z\|_{L^2\to L^2}
		\leq |z|^{\alpha_*-1}$ if $0<|z|\leq1$ and $\|z^{-1}M_z\|_{L^2\to L^2}\leq |z|^{\alpha_*-1}$ if $|z|\geq1$.
		
		For $t\geq1$, choose $\varepsilon=t^{-1}$ in \eqref{eq:vq-contour}. Since $\cos\theta<0$, estimates \eqref{eq:small-resolvent}--\eqref{eq:large-resolvent} give $
		\|v_q(t)\|_{D(A_q)}
		\leq C(\int_{t^{-1}}^1 e^{-c rt}r^{\alpha_*-1}\,dr
		+t^{-\alpha_*} 
		+\int_1^\infty e^{-c rt}
		(1+r^{\alpha^*-\alpha_*})r^{\alpha^*-1}\,dr)
		\leq Ct^{-\alpha_*}.$
		The graph norm of $A_q$ is uniformly equivalent to the $H^2$-norm. Since
		$\partial_t^{\alpha(x)}u=-A_qv_q$, the same estimate gives \eqref{eq:stationary-relaxation}. If $t\geq\tau>0$, the identical contour estimate, now with the exponential factor bounded by $e^{-c\tau r}$ on the high-frequency part, gives \eqref{eq:uniform-terminal-bounds}.
	\end{proof}
	
	\begin{lemma}\label{lema_partial_u}
		Let $q_1,q_2\in\mathcal Q$, and write $u_i=u(\cdot,\cdot;q_i)$. For every $\tau>0$ there is a constant $C_\tau>0$, independent of $q_1$ and $q_2$, such that
		\begin{equation}\label{eq:frac-lipschitz-away-zero}
			\|\partial_t^{\alpha(x)}u_1(t)-\partial_t^{\alpha(x)}u_2(t)\|_{L^2(\Omega)}
			\leq C_\tau\|q_1-q_2\|_{L^2(\Omega)},
			\qquad t\geq\tau.
		\end{equation}
		In addition,
		\begin{equation}\label{eq:frac-lipschitz-decay}
			\|\partial_t^{\alpha(x)}u_1(t)-\partial_t^{\alpha(x)}u_2(t)\|_{L^2(\Omega)}
			\leq Ct^{-\alpha_*}\|q_1-q_2\|_{L^2(\Omega)},
			\qquad t\geq1,
		\end{equation}
		where $C$ is independent of $q_1,q_2$ and $t$.
	\end{lemma}
	\begin{proof}
		Let $U_i$ be the stationary solution associated with $q_i$, put
		$h_i=u_0-U_i$, and set $V=q_1-q_2$. Since $U_1-U_2$ has zero boundary trace and satisfies $A_{q_1}(U_1-U_2)=-VU_2,$ uniform elliptic regularity and the embedding $H^2(\Omega)\hookrightarrow L^\infty(\Omega)$ for $d\leq3$ imply
		\begin{equation}\label{eq:stationary-lipschitz}
			\|U_1-U_2\|_{H^2(\Omega)}+
			\|h_1-h_2\|_{L^2(\Omega)}
			\leq C\|V\|_{L^2(\Omega)}.
		\end{equation}
		
		Let $R_i(z)=(A_{q_i}+M_z)^{-1}$ and $B_i(z)=A_{q_i}R_i(z)=I-M_zR_i(z)$. The resolvent identity reads $R_1(z)-R_2(z)=-R_1(z)V R_2(z).$ 
		For $0<|z|\leq1$, Lemma \ref{lem:frac_deriv_estimate} and
		$H^2(\Omega)\hookrightarrow L^\infty(\Omega)$ give, for $g\in L^2(\Omega)$,
		\begin{align}
			\|(R_1(z)-R_2(z))g\|_{H^2(\Omega)}\leq C\|V R_2(z)g\|_{L^2(\Omega)}\leq C\|V\|_{L^2(\Omega)}\|g\|_{L^2(\Omega)}.\label{eq:small-resolvent-difference}
		\end{align}
		Consequently, since $B_1-B_2=-M_z(R_1-R_2)=M_zR_1VR_2$ and $\|M_z\|_{L^2\to L^2}\le|z|^{\alpha_*}$ for $|z|\le1$, 
		\begin{equation}\label{eq:B-difference-small}
			\|(B_1(z)-B_2(z))g\|_{L^2(\Omega)}
			\leq C|z|^{\alpha_*}\|V\|_{L^2(\Omega)}\|g\|_{L^2(\Omega)}.
		\end{equation}
		For $|z|\geq1$, Proposition~2.1 of \cite{KiaSocYam} and \eqref{eq:large-resolvent} give
		$\|R_i(z)\|_{L^2\to L^2}\leq C|z|^{-\alpha_*}$, and $\|R_i(z)\|_{L^2\to H^2}\leq C|z|^{\alpha^*-\alpha_*}$. Hence the resolvent identity and $H^2(\Omega)\hookrightarrow L^\infty(\Omega)$ imply
		$\|(R_1(z)-R_2(z))g\|_{L^2}\leq C|z|^{\alpha^*-2\alpha_*}\|V\|_{L^2}\|g\|_{L^2}$. Since $B_1-B_2=-M_z(R_1-R_2)$ and $\|M_z\|_{L^2\to L^2}\leq C|z|^{\alpha^*}$, it follows that $\|B_1(z)-B_2(z)\|_{L^2\to L^2}
		\leq C|z|^{2\alpha^*-2\alpha_*}\|V\|_{L^2}$.
		Together with \eqref{eq:stationary-lipschitz} and
		$\|B_2(z)\|_{L^2\to L^2}\leq C(1+|z|^{\alpha^*-\alpha_*})$, this yields
		\begin{equation}\label{eq:large-integrand-difference}
			\|(B_1(z)-B_2(z))z^{-1}M_z h_1
			+B_2(z)z^{-1}M_z(h_1-h_2)\|_{L^2(\Omega)}\leq C|z|^\mu\|V\|_{L^2(\Omega)},
		\end{equation}
		where $\mu=\max\{0,3\alpha^*-2\alpha_*-1\}$.
		
		Because $\partial_t^{\alpha(x)}u_i=-A_{q_i}(u_i-U_i)$, the contour representation gives
		\begin{align*}
			\hspace*{-0.2cm}\partial_t^{\alpha(x)}\big(u_1\!-\!u_2\big) =\!-\frac{1}{2\pi i}\int_{\gamma(\varepsilon,\theta)}e^{zt}
			\Big[(B_1(z)\!-\!B_2(z))z^{-1}M_z h_1
			\!+\!B_2(z)z^{-1}M_z(h_1-h_2)\Big]dz.
		\end{align*}
		By \eqref{eq:B-difference-small}, \eqref{eq:stationary-lipschitz}, and the uniform boundedness of $B_2(z)$ for $|z|\leq1$, the norm of the integrand is bounded by $C|z|^{\alpha_*-1}\|V\|_{L^2(\Omega)}, 0<|z|\leq1.$
		Choosing $\varepsilon=t^{-1}$ for $t\geq1$ and estimating the small-frequency rays and circular arc as in Lemma \ref{lem:frac_deriv_estimate}, while using \eqref{eq:large-integrand-difference} on $|z|\geq1$, yields $\|\partial_t^{\alpha(x)}u_1(t)-\partial_t^{\alpha(x)}u_2(t)\|_{L^2(\Omega)}
		\leq Ct^{-\alpha_*}\|V\|_{L^2(\Omega)}.$
		This proves \eqref{eq:frac-lipschitz-decay}. For $t\geq\tau>0$, the same contour estimate with $e^{-c\tau|z|}$ on the high-frequency part gives \eqref{eq:frac-lipschitz-away-zero}.
	\end{proof}
	
	\section{Uniqueness and conditional stability of inverse problem}\label{section3-IP}
	This section indicates showing the conditional well-posedness of inverse potential problem \eqref{IP-continuous}, including its uniqueness and conditional stability of solutions.
	
	\subsection{Iteration and uniqueness of inverse problem.}
	The reconstruction of the inverse problem is based on the operator $K$, which is defined as 
	\begin{equation}\label{operator}
		K q=\frac{-\partial^{\alpha(x)}_t u(x,T;q)+\Delta \phi(x)+f(x)}{\phi(x)} \;\; {\rm for} \;\; q\in\mathcal{Q}.
	\end{equation}
	For \(q\in\mathcal Q\), Lemma \ref{wellposedness-ivp} gives
	\(\partial_t^{\alpha(x)}u(\cdot,T;q)\in L^2(\Omega)\), and since
	\(\phi\in H^2(\Omega)\), \(f\in L^2(\Omega)\), and \(\phi\ge M_2>0\)
	by Lemma \ref{lemma_positivity}(a), the right-hand side of \eqref{operator} defines an element of
	\(L^2(\Omega)\). Thus \(K\) maps \(\mathcal Q\) into \(L^2(\Omega)\).
	
	For the fixed-point iteration introduced in \eqref{fix-iteration} below we also introduce the working set ${\mathcal D} :=\big\{\psi \in L^2(\Omega): 0\leq \psi(x) \leq \frac{\Delta\phi(x)+f(x)}{\phi(x)} \ \text{a.e. in } \Omega \big\},$
	where the upper bound is understood in the almost-everywhere sense; since \(\phi\ge M_2>0\), \(\Delta\phi\in L^2(\Omega)\) and \(f\in L^2(\Omega)\), the quotient is well defined and belongs to \(L^2(\Omega)\). Lemma \ref{lemma_positivity}(b) applied at \(t=T\) gives $0\le q^*\le q_0:=\frac{\Delta\phi(x)+f(x)}{\phi(x)}
	\quad\text{a.e. in}\Omega,$ so that the exact potential satisfies \(q^*\in\mathcal D\); if moreover \(q_0\in\mathcal Q\), then \(\mathcal D\subseteq\mathcal Q\) and the iterates of \eqref{fix-iteration} remain in \(\mathcal Q\). Note that the uniqueness and stability results of Section \ref{section4} do not rely on the set \(\mathcal D\).
	
	The next lemma contains the equivalence between the fixed point of $K$ and the solution of the inverse problem. 
	\begin{lemma}\label{lemma_equivalence}
		Recalling the equations \eqref{Prob-IBVP} for $u$ and the operator $K$ in \eqref{operator}, the following two statements are equivalent:
		\[
		\textup{(i)}\quad q\in\mathcal Q\text{ is a fixed point of }K;
		\qquad\qquad
		\textup{(ii)}\quad q\in\mathcal Q\text{ satisfies }u(x,T;q)=\phi(x).
		\]
	\end{lemma}
	\begin{proof}
		Assuming $u(x,T;q)=\phi(x)$, it is not difficult to see that $q$ is one fixed point of $K$ and prove the converse. If $q$ is the fixed point, then $q\!=\!\frac{-\partial^{\alpha(x)}_t u(x,T;q)+\Delta \phi(x)+f(x)}{\phi(x)}.$
		Moreover, we note that
		$
		-\partial^{\alpha(x)}_t u(x,T;q)+f(x)=q(x)\phi(x)-\Delta \phi(x)=-\Delta u(x,T;q)+q(x)u(x,T;q).
		$
		Therefore, $v(x)=\phi(x)-u(x,T;q)$ satisfies $-\Delta v(x)+q(x)v(x)=0$ for $x\in \Omega$ with $v(x) =0$ for $x\in\partial\Omega$. By the elliptic maximum principle, \(v \equiv 0\) in \(\Omega\), i.e.  $u(x,T;q)=\phi(x).$ The proof is complete. 
	\end{proof}
	
	The next lemma shows that the operator $K$ is monotone. Note that its proof does not use any commutation between the variable-exponent fractional derivative and the Laplacian; it relies on the positivity of the resolvent family established in Lemma \ref{lemma_positivity}(b).
	\begin{lemma}\label{lemma_monotone}
		Given $q_1,q_2\in {\mathcal Q}$ with $q_1\leq q_2$, we have $Kq_1\le Kq_2$ a.e. in $\Omega$.
	\end{lemma}
	\begin{proof}
		From the definition \eqref{operator} of $K$, $K q_1-K q_2=\frac{\partial^{\alpha(x)}_t u(x,T;q_2)-\partial^{\alpha(x)}_t u(x,T;q_1)}{\phi(x)},$
		so it suffices to prove \(\partial_t^{\alpha(x)}u(\cdot,T;q_1)\ge\partial_t^{\alpha(x)}u(\cdot,T;q_2)\) a.e.\ in \(\Omega\). Set
		\(v=u(\cdot,\cdot;q_1)-u(\cdot,\cdot;q_2)\) and
		\(\widetilde g(t)=(q_2-q_1)u(\cdot,t;q_2)\).
		Subtracting the equations \eqref{Prob-IBVP} for \(q_2\) from those for \(q_1\) shows that \(v\) solves $\partial_t^{\alpha(x)}v+A_{q_1}v=\widetilde g(t)\ge0, v|_{\partial\Omega}=0, v(0)=0,$
		where \(A_{q_1}\) is the Dirichlet realization of \(-\Delta+q_1\); here \(\widetilde g(t)\ge0\) follows from \(q_2\ge q_1\) and Lemma \ref{lemma_positivity}(a).
		
		Let $G_1(s) = (A_{q_1} + M_s)^{-1}$ and let $\nu_{h,\widetilde{h}}$ be the nonnegative measures on $[0,\infty)$ constructed in the proof of Lemma 2.4(b) (with $A_q$ replaced by $A_{q_1}$); recall that the map $h \mapsto \nu_{h,\widetilde{h}}$ is linear and that $\nu_{h,\widetilde{h}} \ge 0$ whenever $h \ge 0$. Denote by $\nu_{h,\widetilde{h}}(\tau) := \nu_{h,\widetilde{h}}([0,\tau])$ the corresponding cumulative distribution function. By the Duhamel formula, for every nonnegative $\widetilde{h} \in L^2(\Omega)$ and every $t \in [0,T]$,
		$\langle v(t), \widetilde{h} \rangle = \int_0^t \nu_{\widetilde{g}(s),\widetilde{h}}(t-s) \, ds \ge 0.$
		
		Moreover, by the representation in the proof of Lemma 2.4(b) (with $G_1$ replaced by $G_2(s) = (A_{q_2} + M_s)^{-1}$ in the construction of $\nu_{h,\widetilde{h}}$), the function $u(\cdot, \cdot; q_2) = u_0 + z_2$, where $z_2$ solves $\partial_t^{\alpha(x)} z_2 + A_{q_2} z_2 = g_2$ with $g_2 = f + \Delta u_0 - q_2 u_0 \ge 0$ and satisfies $\langle z_2(t), \widetilde{h} \rangle = \nu_{g_2,\widetilde{h}}([0,t])$, is nondecreasing in $t$; hence so is $\widetilde{g}(t) = (q_2 - q_1)u(\cdot, t; q_2)$. For $0 \le t \le t'$, split the integral at $t$ to obtain
		$$\langle v(t') - v(t), \widetilde{h} \rangle = \int_0^t \left[ \nu_{\widetilde{g}(s),\widetilde{h}}(t'-s) - \nu_{\widetilde{g}(s),\widetilde{h}}(t-s) \right] ds + \int_t^{t'} \nu_{\widetilde{g}(s),\widetilde{h}}(t'-s) \, ds.$$
		The second integral is nonnegative since $\nu_{\widetilde{g}(s),\widetilde{h}}(\cdot)$ is a nonnegative measure. The first integral is also nonnegative: since $t' - s > t - s$ and the cumulative distribution function $\nu_{h,\widetilde{h}}(\tau)$ is nondecreasing in $\tau$, the integrand is nonnegative. Hence $t \mapsto \langle v(t), \widetilde{h} \rangle$ is nondecreasing in $t$ for all nonnegative $\widetilde{h}$.
		The proof is complete. 
	\end{proof}

	We first record a uniqueness result for ordered fixed points of \(K\), which holds in the variable-exponent setting.
	\begin{lemma}\label{lemma_uniqueness1}
		If $q_1,q_2\in\mathcal Q$ are both fixed points of $K$ with $q_1\le q_2$ a.e. in $\Omega$, then 
		$q_1=q_2$.  
	\end{lemma}
	\begin{proof}
		By Lemma \ref{lemma_equivalence}, we have \(u(x,T;q_1)=u(x,T;q_2)=\phi(x)\). Setting \(v=u(x,t;q_1)-u(x,t;q_2)\), we obtain $\big(\partial_t^{\alpha(x)}-\Delta+q_1\big)v =(q_2-q_1)u(x,t;q_2)\ge0,  (x,t) \in \Omega\times(0,T]$ equipped with homogeneous initial and boundary condition, where the right-hand side is nonnegative by \(q_2\ge q_1\) and Lemma \ref{lemma_positivity}(a). Lemma \ref{lemma_maximum} therefore implies \(v\ge0\) a.e. in \(\Omega\times(0,T]\). Since \(v(\cdot,T)=\phi-\phi=0\), we also have \(\Delta v(\cdot,T)=0\) in \(L^2(\Omega)\).
		
		Evaluating the equation at \(t=T\) gives
		\[
		\partial_t^{\alpha(x)}v(x,T)
		=(q_2-q_1)(x)u(x,T;q_2)
		=(q_2-q_1)(x)\phi(x)\ge0
		\quad\text{a.e. in }\Omega.
		\]
		On the other hand, by the regularity established in Lemma \ref{wellposedness-ivp}, for a.e. \(x\in\Omega\) the function \(v(x,\cdot)\) belongs to \(W^{1,1}(0,T)\) with \(v(x,0)=v(x,T)=0\) (indeed \(v\) is smooth for \(t>0\) and \(\partial_tv(x,t)=\mathcal O(t^{\alpha(x)-1})\) as \(t\to0^+\)). The extremum-type identity for the Caputo derivative \cite[Theorem~1]{yury2009} then yields, for a.e. \(x\in\Omega\),
		\begin{equation*}
			\begin{aligned}
				\partial_t^{\alpha(x)}v(x,T)
				&=\frac{v(x,T)-v(x,0)}{\Gamma(1-\alpha(x))T^{\alpha(x)}}
				+\frac{\alpha(x)}{\Gamma(1-\alpha(x))}
				\int_0^{T}\frac{v(x,T)-v(x,s)}{(T-s)^{1+\alpha(x)}}\,ds \\
				&=-\frac{\alpha(x)}{\Gamma(1-\alpha(x))}
				\int_0^{T}\frac{v(x,s)}{(T-s)^{1+\alpha(x)}}\,ds\le0,
			\end{aligned}
		\end{equation*}
		because \(v\ge0\). Combining the two displays gives $0\!\le (q_2-q_1)(x)\phi(x)=\partial_t^{\alpha(x)}v(x,T)\le0~\text{a.e. in }\Omega,$ hence \((q_2-q_1)\phi=0\) a.e. in \(\Omega\). Since \(\phi\ge M_2>0\) by Lemma \ref{lemma_positivity}(a), we conclude \(q_1=q_2\), which completes the proof.
	\end{proof}

	\begin{lemma}\label{lemma_operator_stability}
		Given $q_1,q_2\in \mathcal Q$, we have $\|Kq_1-Kq_2\|_{L^2(\Omega)}\le C\|q_1-q_2\|_{L^2(\Omega)},$
		where the positive constant $C$ is independent of the choice of $q_1$ and $q_2$.
	\end{lemma}
	\begin{proof}
		From the definition \eqref{operator} of $K$,
		$
		Kq_1-Kq_2=\frac{\partial^{\alpha(x)}_t u(x,T;q_2)-\partial^{\alpha(x)}_t u(x,T;q_1)}{\phi(x)}.
		$
		Since \(\phi=u(\cdot,T;q^*)\ge M_2>0\) by Lemma \ref{lemma_positivity}(a), Lemma \ref{lema_partial_u} gives 
		\begin{align*}
			\|Kq_1- Kq_2\|_{L^2(\Omega)}
			\!\le \!\frac1{M_2} \|\partial^{\alpha(x)}_t u(x,T;q_1)-\partial^{\alpha(x)}_t u(x,T;q_2)\|_{L^2(\Omega)}\!\le\! C\|q_1-q_2\|_{L^2(\Omega)},
		\end{align*}
		which completes the proof. 
	\end{proof}

	Now we state the convergence of the fixed-point iteration; the operator $K$ is monotone by Lemma \ref{lemma_monotone}. The iteration reads
	\begin{equation}\label{fix-iteration}
		q_0=\frac{\Delta\phi(x)+f(x)}{\phi(x)},\ \ q_{n+1}=K q_n, \ n=0,1,2\cdots.
	\end{equation}
	
	\begin{theorem}\label{theorem_uniqueness}
		Let Assumption \ref{assumption} hold, and suppose that \(q_0=(\Delta\phi+f)/\phi\in\mathcal Q\). Then \(\mathcal D\subseteq\mathcal Q\), \(q^*\in\mathcal D\), and the sequence \(\{q_n\}_{n=0}^\infty\) defined by \eqref{fix-iteration} is well defined, decreasing, and converges pointwise a.e. and in \(L^2(\Omega)\) to \(q^*\). 
		
		Consequently, the inverse problem has a unique solution in \(\mathcal Q\); more precisely, if \(q,\tilde q\in\mathcal Q\) are fixed points of \(K\), then \(q=\tilde q\). 
	\end{theorem}
	\begin{proof}
		By Lemma \ref{lemma_positivity}(b) applied at \(t=T\), $q_0-q^*=\frac{\partial_t^{\alpha(x)}u(x,T;q^*)}{\phi(x)}\ge0,$ so that \(0\le q^*\le q_0\); since \(q_0\in\mathcal Q\), this gives \(q^*\in\mathcal D\) and \(\mathcal D\subseteq\mathcal Q\).
		
		We prove by induction that $q^*\le q_{n+1}\le q_n\le q_0, n=0,1,2,\dots.$ For \(n=0\), Lemma \ref{lemma_positivity}(b) applied to \(u(\cdot,\cdot;q_0)\) (which is legitimate because \(q_0\in\mathcal Q\)) gives \(\partial_t^{\alpha(x)}u(x,T;q_0)\ge0\), hence
		$
		q_1=Kq_0
		=\frac{-\partial_t^{\alpha(x)}u(x,T;q_0)+\Delta\phi+f}{\phi}
		=q_0-\frac{\partial_t^{\alpha(x)}u(x,T;q_0)}{\phi}
		\le q_0.
		$
		On the other hand, since \(q^*\le q_0\) and \(K\) is monotone increasing (Lemma \ref{lemma_monotone}), we obtain $q^*=Kq^*\le Kq_0=q_1.$ Hence $q^*\le q_1\le q_0$, which in particular shows \(q_1\in\mathcal D\subseteq\mathcal Q\). The induction step is identical, because \(q_n\in\mathcal Q\) for every \(n\).
		
		Therefore \(\{q_n\}\) is a decreasing sequence bounded below by \(q^*\), hence it converges pointwise a.e. to some \(q_\infty\in L^2(\Omega)\) with \(q^*\le q_\infty\le q_0\), and by the dominated convergence theorem, $\|q_n-q_\infty\|_{L^2(\Omega)}\to0.$ Since \(q_\infty\in\mathcal D\subseteq\mathcal Q\), Lemma \ref{lemma_operator_stability} yields
		\begin{align*}
			\|q_\infty-Kq_\infty\|_{L^2(\Omega)}
			&\le \|q_\infty-q_{n+1}\|_{L^2(\Omega)}
			+\|Kq_n-Kq_\infty\|_{L^2(\Omega)}\\
			&\le \|q_\infty-q_{n+1}\|_{L^2(\Omega)}
			+C\|q_n-q_\infty\|_{L^2(\Omega)}
			\to0
		\end{align*}
		as \(n\to\infty\). Hence \(Kq_\infty=q_\infty\) in \(L^2(\Omega)\), and by Lemma \ref{lemma_equivalence} \(u(\cdot,T;q_\infty)=\phi\), i.e.\ \(q_\infty\) is a solution of the inverse problem. Thus \(q^*\) and \(q_\infty\) are two fixed points of \(K\) in \(\mathcal Q\) with \(q^*\le q_\infty\), and Lemma \ref{lemma_uniqueness1} gives \(q^*=q_\infty\). Therefore the iteration converges to the exact potential \(q^*\).
		
		Finally, let \(\tilde q\in\mathcal Q\) be another solution of the inverse problem. By Lemma \ref{lemma_equivalence}, \(\tilde q\) is a fixed point of \(K\), and Lemma \ref{lemma_positivity}(b) applied at \(t=T\) gives \(\tilde q\le q_0\), so that \(\tilde q\in\mathcal D\). The monotonicity of \(K\) (Lemma \ref{lemma_monotone}) and induction give \(\tilde q=K\tilde q\le Kq_n=q_{n+1}\le q_n\) for every \(n\), so that \(\tilde q\le q_\infty=q^*\). Since \(q^*,\tilde q\in\mathcal Q\) are fixed points of \(K\) with \(\tilde q\le q^*\), Lemma \ref{lemma_uniqueness1} yields \(\tilde q=q^*\). This proves the uniqueness of the fixed point of \(K\) in \(\mathcal Q\), and hence the uniqueness of the solution to the inverse problem.
	\end{proof}
	
	\subsection{Conditional stability of inverse problem}\label{section4}
	
	In this section, we will show the stability of inverse problem \eqref{IP-continuous}, which is concluded as the following theorem.
	
	\begin{theorem}\label{thm-stability}
		Let Assumption \ref{assumption} hold. Then there exists $T_*\geq1$ such that, for every $T\geq T_*$ and every $q_1,q_2\in\mathcal Q$,
		\begin{equation}\label{eq:stability-final}
			\|q_1-q_2\|_{L^2(\Omega)}
			\leq C\|u(\cdot,T;q_1)-u(\cdot,T;q_2)\|_{H^2(\Omega)},
		\end{equation}
		where $C$ is independent of $q_1,q_2$ and $T\geq T_*$.
	\end{theorem}
	\begin{proof}
		For $i=1,2$, the function $z_i=u(\cdot,\cdot;q_i)-u_0$ satisfies
		$
		\partial_t^{\alpha(x)}z_i-\Delta z_i+q_i z_i
		=f+\Delta u_0-q_i u_0\geq(M_1-q_i)u_0\geq0,
		$
		with zero initial and boundary data. Lemma \ref{lemma_maximum} therefore gives
		\begin{equation}\label{eq:uniform-positive-lower-bound}
			u(x,t;q_i)\geq u_0(x)\geq M_2
			\quad\text{for a.e. }(x,t)\in\Omega\times(0,T].
		\end{equation}
		
		For brevity, write $u_i=u(\cdot,T;q_i)$,
		$D_i=\partial_t^{\alpha(x)}u(\cdot,T;q_i)$, and $d=u_1-u_2$. The forward equation gives $q_i=\frac{-D_i+\Delta u_i+f}{u_i},$ and hence $q_1-q_2=J_1+J_2+J_3,$ where $J_1=f\big(\frac1{u_1}-\frac1{u_2}\big),
		J_2=\frac{\Delta u_1}{u_1}-\frac{\Delta u_2}{u_2},
		J_3=\frac{D_2}{u_2}-\frac{D_1}{u_1}.$ The lower bound \eqref{eq:uniform-positive-lower-bound} gives, using that \(f\in C(\overline\Omega)\subset L^\infty(\Omega)\) by Assumption \ref{assumption} and \(|f(1/u_1-1/u_2)|=|f|\,|u_1-u_2|/(u_1u_2)\le M_2^{-2}|f|\,|d|\),
		\begin{equation}\label{eq:J1-estimate}
			\|J_1\|_{L^2(\Omega)}
			\leq M_2^{-2}\|f\|_{L^\infty(\Omega)}\|d\|_{L^2(\Omega)}
			\leq C\|d\|_{L^2(\Omega)}.
		\end{equation}
		For the second term, use $J_2=\frac{\Delta d}{u_1}-\frac{d\,\Delta u_2}{u_1u_2}.$ By Lemma \ref{lem:frac_deriv_estimate}, $\|u_2\|_{H^2(\Omega)}$ is uniformly bounded for $T\!\geq\!1$. Since $H^2(\Omega)\hookrightarrow L^\infty(\Omega)$ for $d\!\leq\!3$, then $\|J_2\|_{L^2(\Omega)}
		\leq C\bigl(\|\Delta d\|_{L^2(\Omega)}
		+\|d\|_{L^\infty(\Omega)}\|\Delta u_2\|_{L^2(\Omega)}\bigr)
		\leq C\|d\|_{H^2(\Omega)}.$
		Finally, $J_3=\frac{D_2-D_1}{u_2}
		+D_1\big(\frac1{u_2}-\frac1{u_1}\big).$ Lemmas \ref{lem:frac_deriv_estimate} and \ref{lema_partial_u} therefore imply
		\begin{align}
			\|J_3\|_{L^2(\Omega)}
			&\leq \frac1{M_2}\|D_2-D_1\|_{L^2(\Omega)}
			+\frac1{M_2^2}\|D_1\|_{L^2(\Omega)}\|d\|_{L^\infty(\Omega)}\nonumber\\
			&\leq CT^{-\alpha_*}\|q_1-q_2\|_{L^2(\Omega)}
			+C\|d\|_{H^2(\Omega)}.\label{eq:J3-estimate}
		\end{align}
		Combining \eqref{eq:J1-estimate}--\eqref{eq:J3-estimate} yields $\|q_1-q_2\|_{L^2(\Omega)}
		\leq C\|d\|_{H^2(\Omega)}
		+CT^{-\alpha_*}\|q_1-q_2\|_{L^2(\Omega)}.$ Choose $T_*\geq1$ so that $CT_*^{-\alpha_*}\leq\frac12$. The last term can then be absorbed into the left-hand side for all $T\geq T_*$, proving \eqref{eq:stability-final}.
	\end{proof}
	
	The stability estimate in Theorem \ref{thm-stability} yields an alternative route to the uniqueness of the inverse problem, independent of the fixed-point iteration and of the condition \(q_0\in\mathcal Q\).
	\begin{corollary}\label{cor-uniqueness}
		Let Assumption \ref{assumption} hold. For every $T\geq T_*$ the inverse problem \eqref{IP-continuous} has at most one solution $q\in\mathcal Q$; equivalently, by Lemma \ref{lemma_equivalence}, the operator $K$ has at most one fixed point in $\mathcal Q$.
	\end{corollary}
	\begin{proof}
		If $q_1,q_2\in\mathcal Q$ satisfy $u(\cdot,T;q_1)=u(\cdot,T;q_2)$, then Theorem \ref{thm-stability} gives $\|q_1-q_2\|_{L^2(\Omega)}=0$, hence $q_1=q_2$.
	\end{proof}
	
	\begin{remark}
		For a variable exponent \(\alpha(x)\), the fractional derivative \(\partial_t^{\alpha(x)}\) and the Laplacian \(\Delta\) do not commute, so the existing proof of Lemma \ref{lemma_monotone} based on the commutation of the two operators does not apply. The proof given here bypasses the commutation issue entirely: the monotonicity of \(K\) follows from the complete monotonicity of the variable-exponent resolvent family \((A_q+s^{\alpha(x)})^{-1}\) established in Lemma \ref{lemma_positivity}(b), which holds without any commutation property. The iteration \eqref{fix-iteration} is the scheme used in Section \ref{Algorithm}.
	\end{remark}

	\section{Error estimates for inverse problem}\label{Algorithm}
	
	In this section, we shall develop a fully discrete scheme for solving the inverse potential problem. To this end, we shall discuss the spatial discretization using the finite element method. A reconstruction algorithm will be presented to recover the potential from the noisy observational data. Finally, we establish the error estimates for the inverse problem \eqref{IP-continuous}. 
	
	\subsection{Finite element scheme for solving the forward problem}
	Let \(\Omega\subset\mathbb R^d\), \(1\le d\le 3\), be a
	bounded Lipschitz domain. The weak formulation of \eqref{Prob-IBVP} reads
	\begin{equation}\label{eq:4.1}
		(\partial_t^{\alpha(x)}u,v)+(\nabla u,\nabla v)+(qu,v)=(f,v),\qquad \forall v\in H_0^1(\Omega),
	\end{equation}
	with \(u(\cdot,t)=u_b\) on \(\partial\Omega\) and \(u(\cdot,0)=u_0\).

	Next, we do the finite element discretization. Let \(\mathcal T_h=\{\tau_e\}_{e=1}^{N_e}\) be a shape-regular simplicial triangulation of \(\Omega\), i.e. $\Omega_h=\bigcup_{e=1}^{N_e}\tau_e$,
	and let $X_h=
	\{v\in C^0(\overline{\Omega}):v|_{\tau_e}\in P_1(\tau_e),\ \forall \tau_e \in \mathcal T_h\}$, and
	$X_h^0=X_h\cap H_0^1(\Omega).$ Here \(P_1(\tau_e)\) denotes the space of linear polynomials on \(\tau_e\). Let
	\(\{P_i\}_{i=1}^{N_0}\) be the nodes of \(\mathcal T_h\) and
	\(\{\phi_i\}_{i=1}^{N_0}\) the corresponding nodal basis functions. Any
	\(u_h\in X_h\) admits the representation $u_h(x,t)=\sum_{i=1}^{N_0}u_i(t)\phi_i(x), u_i(t)=u_h(P_i,t).$
	For \(u_h\in X_h\) and \(v_h\in X_h^0\), the semi-discrete problem is
	\begin{equation}\label{eq:4.2}
		\left(
		\partial_t^{\alpha(x)}u_h(q),v_h
		\right)
		+
		\left(
		\nabla u_h(q),\nabla v_h
		\right)
		+
		\left(
		qu_h(q),v_h
		\right)
		=
		(f,v_h),
	\end{equation}
	with \(u_h(\cdot,t)|_{\partial\Omega}=I_h^\partial u_b\) and
	\(u_h(\cdot,0)=I_hu_0\), where \(I_h\) and \(I_h^\partial\) denote the Lagrange
	interpolation operators in \(X_h\) and on the boundary, respectively.
	
	Define the mass, stiffness, potential matrices and the load vector by
	\[
	M_{ij}=\int_\Omega \phi_i\phi_j\,dx,
	\quad
	K_{ij}=\int_\Omega \nabla\phi_i\cdot\nabla\phi_j\,dx,
	\quad
	Q_{ij}=\int_\Omega q\phi_i\phi_j\,dx,
	\quad
	F_j=\int_\Omega f\phi_j\,dx.
	\]
	Setting $U(t)=\big(u_1(t),\ldots,u_{N_0}(t)\big)^T, \boldsymbol{\alpha}=(\alpha_1,\ldots,\alpha_{N_0})^T,
	\alpha_i=\alpha(P_i),$ and $\partial_t^{\boldsymbol{\alpha}}U(t)
	=
	\big(
	\partial_t^{\alpha_1}u_1(t),\ldots,
	\partial_t^{\alpha_{N_0}}u_{N_0}(t)
	\big)^T,$ we may rewrite \eqref{eq:4.2} as $M\partial_t^{\boldsymbol{\alpha}}U(t)+(K+Q)U(t)=F.$
	
	We next do the time discretization and give fully discrete scheme. Let \(0=t^0<t^1<\cdots<t^N=T\) with \(t^n=n\tau\), \(\tau=T/N\), and write $u_h^n(x)=\sum_{i=1}^{N_0}u_i^n\phi_i(x), U^n=\big(u_1^n,\ldots,u_{N_0}^n\big)^T.$ For each node \(P_i\), set \(\alpha_i=\alpha(P_i)\). The Caputo derivative at
	\(P_i\) is approximated by $\bar{\partial}_{\tau}^{\alpha_i}u_i^n=
	\frac{1}{\Gamma(2-\alpha_i)\tau^{\alpha_i}}
	\big[u_i^n+\sum_{j=0}^{n-2}
	\big(d_{j+1}^{(i)}-d_j^{(i)}\big)u_i^{\,n-1-j}-d_{n-1}^{(i)}u_i^0\big],$
	where $d_j^{(i)}=(j+1)^{1-\alpha_i}-j^{1-\alpha_i}, j=0,1,\ldots,n-1.$ Let
	$
	\lambda_i=\frac{1}{\Gamma(2-\alpha_i)\tau^{\alpha_i}},
	$ and
	$
	H_i^n=
	d_{n-1}^{(i)}u_i^0
	-
	\sum_{j=0}^{n-2}
	\big(d_{j+1}^{(i)}-d_j^{(i)}\big)
	u_i^{\,n-1-j},
	$
	so that $\bar{\partial}_{\tau}^{\alpha_i}u_i^n=\lambda_i(u_i^n-H_i^n).$ Define $\lambda_h(x)=\sum_{i=1}^{N_0}\lambda_i\phi_i(x),
	H_h^n(x)=\sum_{i=1}^{N_0}H_i^n\phi_i(x)$, $\alpha_h(x)=\sum_{i=1}^{N_0}\alpha_i\phi_i(x).$ The fully discrete scheme reads: for \(n=1,\ldots,N\), find \(u_h^n\in X_h\)
	such that, for all \(v_h\in X_h^0\),
	\begin{equation}\label{eq:4.4}
		(\bar{\partial}_\tau^{\alpha_h}u_h^n(q),v_h)
		+
		(\nabla u_h^n(q),\nabla v_h)
		+
		(qu_h^n(q),v_h)
		=
		(f,v_h),
	\end{equation}
	or equivalently
	\begin{equation}\label{eq:4.5}
		(\lambda_hu_h^n(q),v_h)
		+
		(\nabla u_h^n(q),\nabla v_h)
		+
		(qu_h^n(q),v_h)
		=
		(f,v_h)
		+
		(\lambda_hH_h^n(q),v_h),
	\end{equation}
	with \(u_h^0=I_hu_0\) and \(\gamma_0(u_h^n)=I_h^\partial u_b\). Introducing $M_{\alpha,ij}=\int_\Omega \lambda_h\phi_i\phi_j\,dx,$ the matrix form of \eqref{eq:4.5} is $(M_\alpha+K+Q)U^n=F+M_\alpha H^n,$ where $U^n=\big(u_1^n, u_2^n,\ldots, u_{N_0}^n\big)^T, H^n=\big(H_1^n,H_2^n,\ldots,H_{N_0}^n\big)^T.$
	
	\subsection{Error estimate for the forward problem}
	
	Define the bilinear form $A(u,v)=\left(\nabla u,\nabla v\right)+(qu,v).$ Assume $0\le q(x)\le M_1$, so that \(A\) is uniformly elliptic on \(H_0^1(\Omega)\): $A(v,v)\ge c\|\nabla v\|_{L^2(\Omega)}^2.$
	Let \(I_h\) and \(R_h\) denote the Lagrange interpolation and the Ritz projection,
	respectively. For \(u\in H^2(\Omega)\), we have the esimates
	\begin{align}
		&\|I_hu-u\|_{L^2(\Omega)}
		+ h\|\nabla(I_hu-u)\|_{L^2(\Omega)}
		\le Ch^2\|u\|_{H^2(\Omega)}, \label{eq:4.7}\\[4pt]
		&\|u-R_hu\|_{L^2(\Omega)}
		\le Ch^2\|u\|_{H^2(\Omega)}, \quad
		\|\partial_t(u-R_hu)\|_{L^2(\Omega)}
		\le Ch^2\|\partial_tu\|_{H^2(\Omega)}. \label{eq:4.8}
	\end{align}
	
	\begin{theorem}\label{thm:forward-error}
		Let \(q^*\) be the exact potential and \(u=u(q^*)\) the exact solution
		of \eqref{Prob-IBVP}. Under the above discretization and regularity
		assumptions, for any \(0<\varepsilon<1\) there exists a constant \(C>0\),
		independent of \(h\) and \(\tau\), such that
		\begin{equation}\label{eq:4.9}
			\|u^n-u_h^n(q^*)\|_{L^2(\Omega)}
			\le C\Bigl[ h^2+h^{2-\varepsilon}\max_{x\in\Omega}\big\{t_n^{-\alpha(x)},t_n^{-(1-\varepsilon)\alpha(x)}\big\}
			+\tau\max_{x\in\Omega}t_n^{-\alpha(x)-1}\Bigr]
		\end{equation}
		for $n=1,\ldots,N$.
	\end{theorem}
	
	\begin{proof}
		Write \(u_h^n=u_h^n(q^*)\), \(A=A(1,q^*)\), and split $u^n-u_h^n
		=
		(u^n-R_hu^n)+(R_hu^n-u_h^n)
		=:
		\rho^n+\theta^n.$
		By \eqref{eq:4.8},  $\|\rho^n\|_{L^2(\Omega)}\le Ch^2\|u^n\|_{H^2(\Omega)}\le Ch^2.$ Taking \(t=t_n\) in the continuous weak formulation and using the definition of
		\(R_h\), we obtain
		\begin{equation}\label{eq:4.10}
			(\partial_t^{\alpha(x)}u^n,v_h)+A(R_hu^n,v_h)=(f,v_h),
			\qquad \forall v_h\in X_h^0.
		\end{equation}
		Subtracting \eqref{eq:4.10} from \eqref{eq:4.4} and substituting
		\(u_h^n=R_hu^n-\theta^n\) gives
		\begin{equation}\label{eq:4.11}
			(\bar{\partial}_{\tau}^{\alpha_h}\theta^n,v_h)
			+
			A(\theta^n,v_h)
			=
			(\omega^n+\eta^n,v_h),
			\qquad \forall v_h\in X_h^0,
		\end{equation}
		where $\omega^n=\bar{\partial}_{\tau}^{\alpha_h}R_hu^n-\bar{\partial}_{\tau}^{\alpha_h}u^n, \eta^n
		=\bar{\partial}_{\tau}^{\alpha_h}u^n-\partial_t^{\alpha(x)}u^n.$ Recall that $\bar{\partial}_{\tau}^{\alpha_h}\theta^n=\lambda_h(\theta^n-H_h^n(\theta)),$ where \(H_h^n(\theta)\) is the history term determined by
		\(\theta^0,\theta^1,\dots,\theta^{n-1}\). Taking \(v_h=\theta^n\) in
		\eqref{eq:4.11} yields $(\lambda_h\theta^n,\theta^n)+A(\theta^n,\theta^n)=(\lambda_hH_h^n(\theta),\theta^n)+(\omega^n+\eta^n,\theta^n).$ Since \(\lambda_h\ge c\tau^{-\alpha_*}\) and \(A(\theta^n,\theta^n)\ge
		c\|\nabla\theta^n\|_{L^2}^2\), we obtain $\|\theta^n\|_{L^2(\Omega)}
		\le
		C\big(
		\|\theta^0\|_{L^2(\Omega)}
		+
		\sum_{j=1}^n\tau(t_n-t_j)^{-\alpha_*}
		\|\omega^j+\eta^j\|_{L^2(\Omega)}
		\big),$
		where we have used the standard discrete fractional Grönwall inequality
		for the variable-exponent Caputo derivative; see, e.g.,
		\cite[Lemma 3.1]{JinMC} or \cite[Lemma 2.2]{JinNM}. The two consistency
		terms satisfy the key estimates $ 
		\|\eta^j\|_{L^2(\Omega)}
		\le
		C\tau
		\max_{x\in\Omega}
		t_j^{-\alpha(x)-1},$ and $\|\omega^j\|_{L^2(\Omega)}
		\le Ch^{2-\varepsilon}
		\max_{x\in\Omega}
		\{
		t_j^{-\alpha(x)},
		t_j^{-(1-\varepsilon)\alpha(x)}
		\}$.
		
		Moreover, \(u_h^0=I_hu_0\) and hence $\|\theta^0\|_{L^2(\Omega)}
		=
		\|R_hu_0-I_hu_0\|_{L^2(\Omega)}
		\le
		Ch^2,$ by \eqref{eq:4.7}--\eqref{eq:4.8}. Substituting these bounds into the
		Grönwall inequality and summing up the discrete convolution, we arrive at
		\[
		\|\theta^n\|_{L^2(\Omega)}
		\le
		C\Bigl[
		h^2+
		h^{2-\varepsilon}
		\max_{x\in\Omega}
		\big\{
		t_n^{-\alpha(x)},
		t_n^{-(1-\varepsilon)\alpha(x)}
		\big\}+\tau\max_{x\in\Omega}
		t_n^{-\alpha(x)-1}
		\Bigr].
		\]
		Combining this with the estimate for \(\rho^n\) and applying the triangle
		inequality yields \eqref{eq:4.9}.
	\end{proof}

	%\begin{remark}
	%The proof isolates the two error mechanisms: the spatial error is governed by
	%the Ritz projection estimate \(O(h^2)\), while the temporal error is governed by
	%the discrete convolution weights of the variable-exponent Caputo derivative.
	%The non-commutativity of \(\partial_t^{\alpha(x)}\) and \(\Delta\) only affects the
	%constants and the exponents \(\alpha(x)\), not the structure of the estimate.
	%\end{remark}
	
	\subsection{Error estimate for the inverse problem}\label{section4.3}
	
	In this subsection, we follow the idea of \cite[Lemma~4.10]{zhouzhi2022} to construct a stable finite element approximation \(\psi_h\) of \(\Delta\phi\) from the noisy terminal observation \(\phi^\delta\), and then directly
	establish the error estimate for the inverse problem.
	
	We assume that the noisy terminal data \(\phi^\delta\) satisfies
	\begin{equation}\label{eq:noise-C}
		\phi^\delta\in C(\overline\Omega),\qquad
		\phi^\delta\ge M_2,\qquad
		\gamma_0(\phi^\delta)=\gamma_0(\phi)=u_b,
		\qquad
		\|\phi^\delta-\phi\|_{C(\overline\Omega)}=\delta.
	\end{equation}
	Moreover, we assume that \(q^*|_{\partial\Omega}\) is known, so that
	\(\gamma_0(\Delta\phi)=\gamma_0(q^*\phi-f)
	=\gamma_0(q^*)u_b-\gamma_0(f)\) is available.
	
	Let \(X_h\) be the finite element space introduced in
	Section~\ref{Algorithm}, and let \(I_h^0\) denote the boundary interpolation
	operator. Define \(\psi_h\in X_h\) by
	\begin{equation}\label{eq:psi-h}
		\gamma_0(\psi_h)
		=
		I_h^0\bigl(\gamma_0(q^*)u_b-\gamma_0(f)\bigr),
		\quad (\psi_h,v_h)=-(\nabla I_h\phi^\delta,\nabla v_h),
		\quad
		\forall v_h\in X_h^0.
	\end{equation}
	%This is the direct analogue of \cite[eq.~(4.21)]{zhouzhi2022}.
	
	\begin{lemma}\label{lem:psi-h}
		Let Assumption \ref{assumption} hold and let
		\eqref{eq:noise-C} be valid. Then there exists a constant \(C>0\),
		independent of \(h\) and \(\delta\), such that
		\begin{equation}\label{eq:psi-h-error}
			\|\psi_h-\Delta\phi\|_{L^2(\Omega)}
			\le
			C\big(\frac{\delta}{h^2}+h\big).
		\end{equation}
	\end{lemma}
	
	\begin{proof}
		The proof is the same as that of \cite[Lemma~4.10]{zhouzhi2022}.
		For completeness, we sketch it. Let \(\widetilde\psi_h\in X_h\) be defined by
		\eqref{eq:psi-h} with \(\phi^\delta\) replaced by the exact data \(\phi\).
		Then \(\psi_h-\widetilde\psi_h\in X_h^0\), and $(\psi_h-\widetilde\psi_h,v_h)
		=
		-(\nabla I_h(\phi^\delta-\phi),\nabla v_h), \forall v_h\in X_h^0.$
		Using the inverse inequality and 
		\(\|I_h(\phi^\delta-\phi)\|_{L^2(\Omega)}\le C\delta\), we obtain $\|\psi_h-\widetilde\psi_h\|_{L^2(\Omega)}
		\le
		C h^{-2}\|I_h(\phi^\delta-\phi)\|_{L^2(\Omega)}
		\le
		C\frac{\delta}{h^2}.$
		On the other hand, by the finite element approximation property, the
		superconvergence estimate, and \(\Delta\phi\in W^{1,p}(\Omega)\) for some
		\(p>\max(d,2)\), the same argument as in \cite[Lemma~4.10]{zhouzhi2022}
		gives $\|\widetilde\psi_h-\Delta\phi\|_{L^2(\Omega)}
		\le
		C h.$ Combining the above two estimates yields \eqref{eq:psi-h-error}.
	\end{proof}
	
	Now we define the approximate fixed-point operator
	\(K^\delta_{h,\tau}:\mathcal Q\to\mathcal Q\) by
	\begin{equation}\label{operator-noise}
		K^\delta_{h,\tau}q
		=
		\Pi_{[0,M_1]}
		\left(
		\frac{f(x)+\psi_h(x)
			-\bar{\partial}_{\tau}^{\alpha_h}u_h^N(x;q)}
		{\phi^\delta(x)}
		\right),
	\end{equation}
	where \(\Pi_{[0,M_1]}\) denotes the pointwise projection onto \([0,M_1]\).
	The corresponding iteration reads
	\begin{equation}\label{fix-iteration-step2}
		q^\delta_0
		=
		\Pi_{[0,M_1]}
		\left(
		\frac{f(x)+\psi_h(x)}{\phi^\delta(x)}
		\right),
		\qquad
		q_{n+1}^\delta=K^\delta_{h,\tau}q_n^\delta,
		\quad n=0,1,2,\dots.
	\end{equation}
	
	\begin{lemma}\label{lem:contraction-Kdelta}
		Let Assumption \ref{assumption} hold, and let \eqref{eq:noise-C} and
		\eqref{eq:psi-h-error} be satisfied with \(\delta<M_2/2\). Then there exists
		a constant \(C>0\), independent of \(h,\tau,\delta\), such that for all
		\(q_1,q_2\in\mathcal Q\),
		\begin{equation}\label{eq:contraction-Kdelta}
			\|K^\delta_{h,\tau}q_1-K^\delta_{h,\tau}q_2\|_{L^2(\Omega)}
			\le
			C\max\{T^{-\alpha_*},T^{-(1-\epsilon)\alpha_*}\}
			\|q_1-q_2\|_{L^2(\Omega)}
		\end{equation}
		for any \(\epsilon\in(0,\min(1,2-d/2))\).
		In particular, if \(T\ge T_*\) is sufficiently large so that $C\max\{T_*^{-\alpha_*},T_*^{-(1-\epsilon)\alpha_*}\}\le\kappa<1,$ then \(K^\delta_{h,\tau}\) is a contraction on \(\mathcal Q\subset L^2(\Omega)\) with contraction constant \(\kappa\). Consequently, the iteration
		\eqref{fix-iteration-step2} has a unique fixed point
		\(q^\delta\in\mathcal Q\), and $\|q^\delta_{n+1}-q^\delta\|_{L^2(\Omega)}
		\le
		\kappa\|q^\delta_n-q^\delta\|_{L^2(\Omega)}$, for $n=0,1,2,\dots.$
	\end{lemma}
	
	\begin{proof}
		By \eqref{eq:noise-C}, \(\phi^\delta\ge M_2\). Since \(\Pi_{[0,M_1]}\) is
		nonexpansive and \(\phi^\delta\ge M_2\), we have $|K^\delta_{h,\tau}q_1-K^\delta_{h,\tau}q_2|
		\le
		\frac{1}{M_2}
		\left|
		\bar{\partial}_{\tau}^{\alpha_h}
		\bigl(u_h^N(q_1)-u_h^N(q_2)\bigr)
		\right|.$ Taking the \(L^2(\Omega)\)-norm and using the fully discrete analogue of
		Lemma \ref{lema_partial_u} gives $\left\|
		\bar{\partial}_{\tau}^{\alpha_h}
		\bigl(u_h^N(q_1)-u_h^N(q_2)\bigr)
		\right\|_{L^2(\Omega)}
		\le
		C\max\{T^{-\alpha_*},T^{-(1-\epsilon)\alpha_*}\}
		\|q_1-q_2\|_{L^2(\Omega)}.$ Combining the two inequalities yields \eqref{eq:contraction-Kdelta}.
		Banach's fixed point theorem then gives the unique fixed point and the
		linear convergence.
	\end{proof}
	
	Now we present the main error estimate for the inverse problem.
	
	\begin{theorem}\label{thm:inverse-error}
		Let Assumption \ref{assumption} hold, and let the noise assumption
		\eqref{eq:noise-C} be satisfied. Let \(q^\delta\in\mathcal Q\) be a fixed
		point of the discrete operator \(K^\delta_{h,\tau}\) defined in
		\eqref{operator-noise}. Then there exists \(T_*>0\) such that for every
		\(T\ge T_*\),
		\begin{equation}\label{eq:inverse-error-general}
			\|q^*-q^\delta\|_{L^2(\Omega)}
			\le
			C\left(\frac{\delta}{h^2}+h+\tau\right),
		\end{equation}
		where \(C>0\) is independent of \(h,\tau,\delta\) and \(T\ge T_*\).
		In particular, if the discretization parameters are chosen as $h\sim \tau\sim \delta^{1/3}$,
		then the reconstruction error satisfies $\|q^*-q^\delta\|_{L^2(\Omega)}
		=
		\mathcal O(\delta^{1/3}).$
	\end{theorem}
	
	\begin{proof}
		Fixed-point equations for \(q^*\) and \(q^\delta\) yeild $q^*
		=
		\frac{f+\Delta\phi-\partial_t^{\alpha(x)}u(\cdot,T;q^*)}{\phi},$ and $q^\delta
		=
		\Pi_{[0,M_1]}
		\left(
		\frac{f+\psi_h-\bar{\partial}_{\tau}^{\alpha_h}u_h^N(\cdot;q^\delta)}
		{\phi^\delta}
		\right).$ Since \(\Pi_{[0,M_1]}\) is nonexpansive, we obtain
		\[
		\|q^*-q^\delta\|_{L^2(\Omega)}
		\le
		\big\|
		\frac{f+\Delta\phi-\partial_t^{\alpha(x)}u(\cdot,T;q^*)}{\phi}
		-
		\frac{f+\psi_h-\bar{\partial}_{\tau}^{\alpha_h}u_h^N(\cdot;q^\delta)}
		{\phi^\delta}
		\big\|_{L^2(\Omega)}.
		\]
		We split the right-hand side into two parts: $\|q^*-q^\delta\|_{L^2(\Omega)}
		\le I_1+I_2,$ where $I_1
		=
		\big\|
		\frac{f+\Delta\phi}{\phi}
		-
		\frac{f+\psi_h}{\phi^\delta}
		\big\|_{L^2(\Omega)},$ and $I_2
		=
		\big\|
		\frac{\partial_t^{\alpha(x)}u(\cdot,T;q^*)}{\phi}
		-
		\frac{\bar{\partial}_{\tau}^{\alpha_h}u_h^N(\cdot;q^\delta)}
		{\phi^\delta}
		\big\|_{L^2(\Omega)}.$
		
		For \(I_1\), using \eqref{eq:noise-C}, \eqref{eq:psi-h-error},
		\(\phi\ge M_2\), \(\phi^\delta\ge M_2\), and
		\(f+\Delta\phi\in L^2(\Omega)\), we get $I_1
		\le
		C\big(\|\Delta\phi-\psi_h\|_{L^2(\Omega)}
		+\|\phi-\phi^\delta\|_{C(\overline\Omega)}\big)
		\le
		C\left(\frac{\delta}{h^2}+h\right).$
		For \(I_2\), we first note that \(\phi^\delta\ge M_2\). Hence $I_2
		\le
		C\big\|
		\bar{\partial}_{\tau}^{\alpha_h}u_h^N(\cdot;q^\delta)
		-\partial_t^{\alpha(x)}u(\cdot,T;q^*)
		\big\|_{L^2(\Omega)}
		+
		C\delta.$ Furthermore\[\begin{aligned}  &\big\|\bar{\partial}_{\tau}^{\alpha_h}u_h^N(\cdot;q^\delta)
			-\partial_t^{\alpha(x)}u(\cdot,T;q^*)
			\big\|_{L^2(\Omega)} \le\left\|\bar{\partial}_{\tau}^{\alpha_h}u_h^N(\cdot;q^\delta)-\bar{\partial}_{\tau}^{\alpha_h}u_h^N(\cdot;q^*)
			\right\|_{L^2(\Omega)} \\	&+\left\|\bar{\partial}_{\tau}^{\alpha_h}u_h^N(\cdot;q^*)-\bar{\partial}_{\tau}^{\alpha_h}u^N(\cdot;q^*)\right\|_{L^2(\Omega)}+\big\|\bar{\partial}_{\tau}^{\alpha_h}u^N(\cdot;q^*)
			-\partial_t^{\alpha(x)}u(\cdot,T;q^*)
			\big\|_{L^2(\Omega)}.
		\end{aligned}
		\]
		The first term is controlled by the discrete analogue of
		Lemma \ref{lema_partial_u}: $\|
		\bar{\partial}_{\tau}^{\alpha_h}u_h^N(\cdot;q^\delta)
		-\bar{\partial}_{\tau}^{\alpha_h}u_h^N(\cdot;q^*)
		\|_{L^2(\Omega)}
		\le
		C\max\{T^{-\alpha_*},T^{-(1-\epsilon)\alpha_*}\}
		\|q^\delta-q^*\|_{L^2(\Omega)}.$ The second term is the forward spatial discretization error, which by
		Theorem \ref{thm:forward-error} is bounded by \(Ch\). The third term is
		the time-discretization error of the variable-exponent Caputo derivative,
		which is bounded by \(C\tau\). Hence $I_2
		\le
		C\left(\frac{\delta}{h^2}+h+\tau\right)
		+
		C\max\{T^{-\alpha_*},T^{-(1-\epsilon)\alpha_*}\}
		\|q^*-q^\delta\|_{L^2(\Omega)}.$ Combining the estimates for \(I_1\) and \(I_2\), we arrive at
		$\|q^*-q^\delta\|_{L^2(\Omega)}
		\le
		C\left(\frac{\delta}{h^2}+h+\tau\right)
		+
		C\max\{T^{-\alpha_*},T^{-(1-\epsilon)\alpha_*}\}
		\|q^*-q^\delta\|_{L^2(\Omega)}.$
		Choose \(T_*\ge1\) large enough so that $C\max\{T_*^{-\alpha_*},T_*^{-(1-\epsilon)\alpha_*}\}\le\frac12.$ Then the last term can be absorbed into the left-hand side for all
		\(T\ge T_*\), which proves \eqref{eq:inverse-error-general}.
		
		Finally, taking \(h\sim\tau\sim\delta^{1/3}\) gives
		$\frac{\delta}{h^2}\sim\delta^{1/3},\; h\sim\delta^{1/3},\;
		\tau\sim\delta^{1/3}$, and therefore $\|q^*-q^\delta\|_{L^2(\Omega)}
		=
		\mathcal O(\delta^{1/3}).$ The proof is complete.
	\end{proof}

	\section{Numerical reconstruction for the inverse problem}\label{sec5}
	
	We consider the domain \(\Omega=(0,4)^2\) and the terminal time \(T=5\). The problem data are chosen as $f(x,y)=30, u_b(x,y)=\frac{x(4-x)}{5}+3, u_0(x,y)=x(4-x)\big(\frac14+\frac{y(4-y)}{40}\big)+3.$ These data satisfy Assumption~\ref{assumption}. To examine the influence of the variable exponent, we test three cases: $\alpha_1(x,y)=0.5,\alpha_2(x,y)=0.5+\frac{xy}{40}, \alpha_3(x,y)=0.5+0.25\sin(\pi x)\sin(\pi y).$ Thus \(\alpha(x,y)\) varies in space and remains strictly between \(0\) and \(1\), so that the variable-exponent Caputo derivative is well defined.
	
	The noisy terminal observation is generated by $\phi^{\delta}(x_i)=\phi(x_i)+\delta\,\zeta(x_i)$, where $\zeta$ follows the standard uniform distribution in $[-1,1]$, and \(x_i\) are grid points of a fine partition of \(\Omega\). Then, to compute the numerical reconstruction, we follow the idea in Section \ref{Algorithm} and design the Algorithm \ref{alg:recon} below. 
	
	\begin{algorithm}[H]
		\caption{Iterative algorithm for finding \(q^\ast(x)\)}\label{alg:recon}
		\begin{algorithmic}
			\REQUIRE Exponent \(\alpha(x)\), terminal time \(T\), source term
			\(f(x)\), initial condition \(u_0(x)\), boundary data \(u_b(x)\), noisy
			observational data \(\phi^\delta\), upper bound $M_1$ of the potential and discretization \(h\) and \(\tau\).
			\ENSURE Approximate potential \(q^\ast\).
			
			\STATE Compute \(\psi_h\) by \eqref{eq:psi-h}; set
			\[
			q_0
			=
			\Pi_{[0,M_1]}
			\left(
			\frac{f(x)+\psi_h(x)}{\phi^\delta(x)}
			\right),
			\quad k=0.
			\]
			\WHILE{\(\|q_{k+1}-q_k\|_{L^2(\Omega)}/\|q_k\|_{L^2(\Omega)} > 10^{-10}\)}
			\STATE Compute \(u_h^N(\cdot;q_k)\) by the fully discrete scheme
			\eqref{eq:4.4}; update
			\[
			q_{k+1}
			=
			K_{h,\tau}^\delta q_k
			=
			\Pi_{[0,M_1]}
			\left(
			\frac{f(x)+\psi_h(x)
				-\bar{\partial}_{\tau}^{\alpha_h}u_h^N(x;q_k)}
			{\phi^\delta(x)}
			\right).
			\]
			\ENDWHILE
			
			\STATE  \(q^\ast \leftarrow q_{k+1}\).
			
			\ENSURE Approximate potential \(q^\ast\).
		\end{algorithmic}
	\end{algorithm}

\begin{figure}[H]
		\centering
		\begin{subfigure}[t]{0.31\textwidth}
			\centering
			\includegraphics[width=\linewidth]{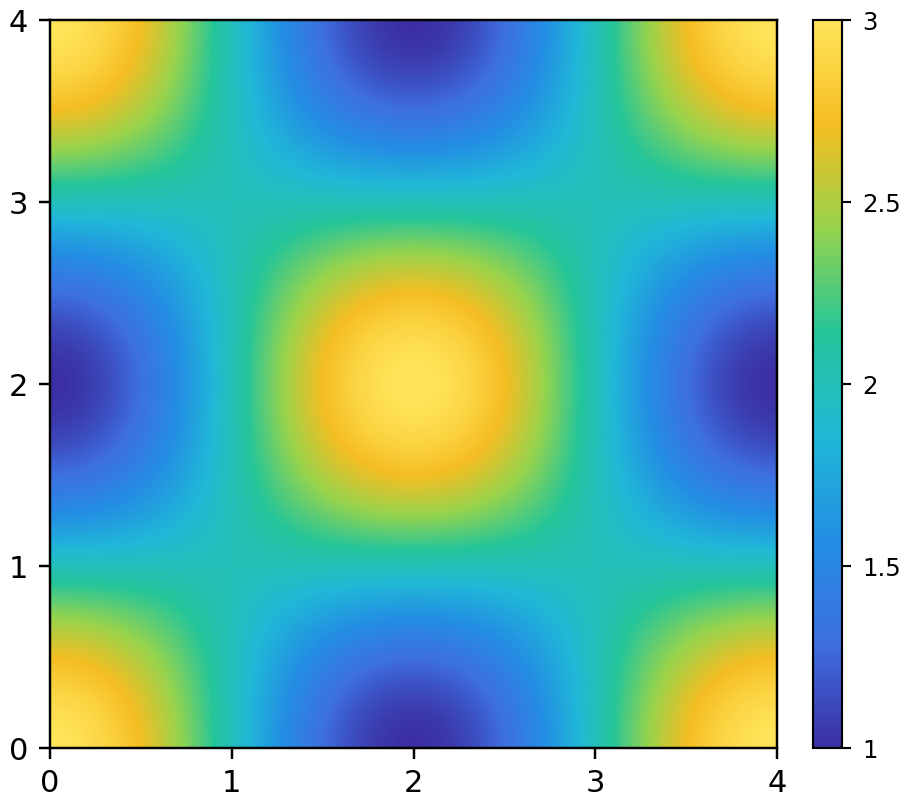}
			\caption{\(q_1^*\)}
		\end{subfigure}
		\hspace{0.01\textwidth}
		\begin{subfigure}[t]{0.31\textwidth}
			\centering
			\includegraphics[width=\linewidth]{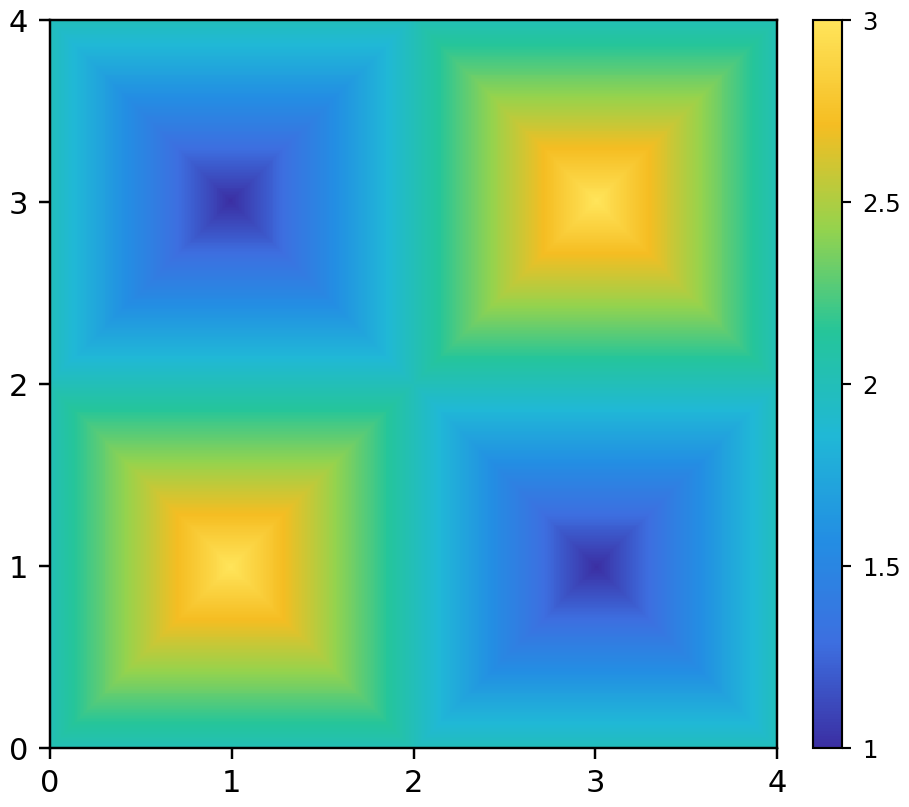}
			\caption{\(q_2^*\)}
		\end{subfigure}
		\hspace{0.01\textwidth}
		\begin{subfigure}[t]{0.31\textwidth}
			\centering
			\includegraphics[width=\linewidth]{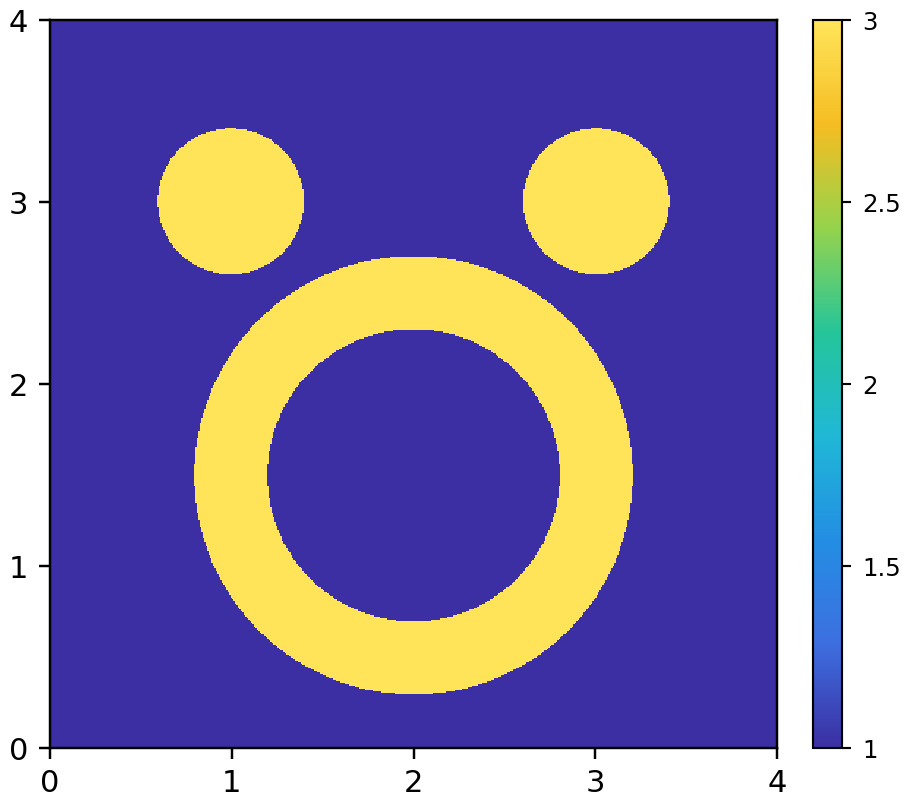}
			\caption{\(q_3^*\)}
		\end{subfigure}
		\caption{Profiles of exact solutions \(q_1^*\), \(q_2^*\) and \(q_3^*\) in Example~\ref{ex1}.}
		\label{ex1-fig1}
	\end{figure}

	According to Theorem~\ref{thm:inverse-error}, the reconstruction error satisfies $\|q^*-q^\delta\|_{L^2(\Omega)}
	\le
	C\left(\frac{\delta}{h^2}+h+\tau\right),$ where \(q^*\) is the exact potential and \(q^\delta\) is its reconstruction by Algorithm \ref{alg:recon}. Hence, if the discretization parameters are chosen as $h\sim \tau\sim \delta^{1/3}$, then the expected convergence rate is $\|q^*-q^\delta\|_{L^2(\Omega)}
	=
	\mathcal O(\delta^{1/3}).$ In the experiments below, we therefore take $h=\delta^{1/3}, \tau=\frac{h}{15}$, and measure the relative error $Err:=
	\frac{\|q^\delta-q^*\|_{L^2(\Omega)}}{\|q^*\|_{L^2(\Omega)}}.$
	\begin{example}\label{ex1}
		We test three exact potentials with different regularity properties: a smooth potential $q_1^*$, a piecewise smooth potential $q_2^*$, and a discontinuous potential $q_3^*$.
		Their profiles are shown in Figure~\ref{ex1-fig1}. The first one is smooth, the second is continuous but has piecewise smooth derivatives, and the third has jump discontinuities. These examples allow us to examine the effect of the regularity of \(q^*\) on the reconstruction.
	\end{example}
	%\begin{itemize}
	%\item[(1)] a smooth potential \(q_1^*\);
	%\item[(2)] a piecewise smooth potential \(q_2^*\);
	%\item[(3)] a discontinuous potential \(q_3^*\).
	%\end{itemize}
	%========================= Figure 1 =========================
	
	The numerical results for the variable exponent \(\alpha=\alpha_2\) are reported in Table~\ref{ex1-tab1}. For the smooth and piecewise smooth potentials, the reconstruction error decreases as the noise level \(\delta\) decreases. For example, for \(q_1^*\), the error decreases from \(9.77\times 10^{-2}\) at \(1\%\) noise to \(2.90\times 10^{-2}\) at \(0.1\%\) noise and further to \(1.91\times 10^{-2}\) at \(0.01\%\) noise. A similar behavior is observed for \(q_2^*\). This is consistent with the stability estimate in Theorem~\ref{thm-stability} and the error estimate in Theorem~\ref{thm:inverse-error}: a smaller noise leads to a more accurate reconstruction, provided that \(T\) is sufficiently large.
	
	For the discontinuous potential \(q_3^*\), the error is significantly larger and decreases more slowly. At \(1\%\), \(0.1\%\), and \(0.01\%\) noise, the errors are \(1.69\times 10^{-1}\), \(1.10\times 10^{-1}\), and \(9.15\times 10^{-2}\), respectively. This behavior is expected, since \(q_3^*\) does not satisfy the continuity assumption in Assumption~\ref{assumption}(a). The result illustrates the necessity of the regularity assumption on the exact potential in the theoretical analysis.

	%............................................table 1...........................................
	\begin{table}[h!]
		\caption{Reconstruction error \(Err\) for Example~\ref{ex1} with \(\alpha=\alpha_2\).}\label{ex1-tab1}
		\centering
		\renewcommand{\arraystretch}{0.95}   
		\setlength{\tabcolsep}{4pt}         
		\begin{tabular}{cccc}
			\toprule
			Exact potential & \(\delta=1\%\) & \(\delta=0.1\%\) & \(\delta=0.01\%\) \\
			\midrule
			Smooth potential \(q_1^*\)          & 9.77e-2 & 2.90e-2 & 1.91e-2 \\
			Piecewise smooth potential \(q_2^*\) & 1.05e-1 & 3.22e-2 & 1.94e-2 \\
			Discontinuous potential \(q_3^*\)   & 1.69e-1 & 1.10e-1 & 9.15e-2 \\
			\bottomrule
		\end{tabular}
	\end{table}
	
\begin{figure}
		\centering
		\begin{subfigure}[t]{0.23\textwidth}
			\centering
			\caption*{(a) noise 1\%}
			\includegraphics[width=\linewidth]{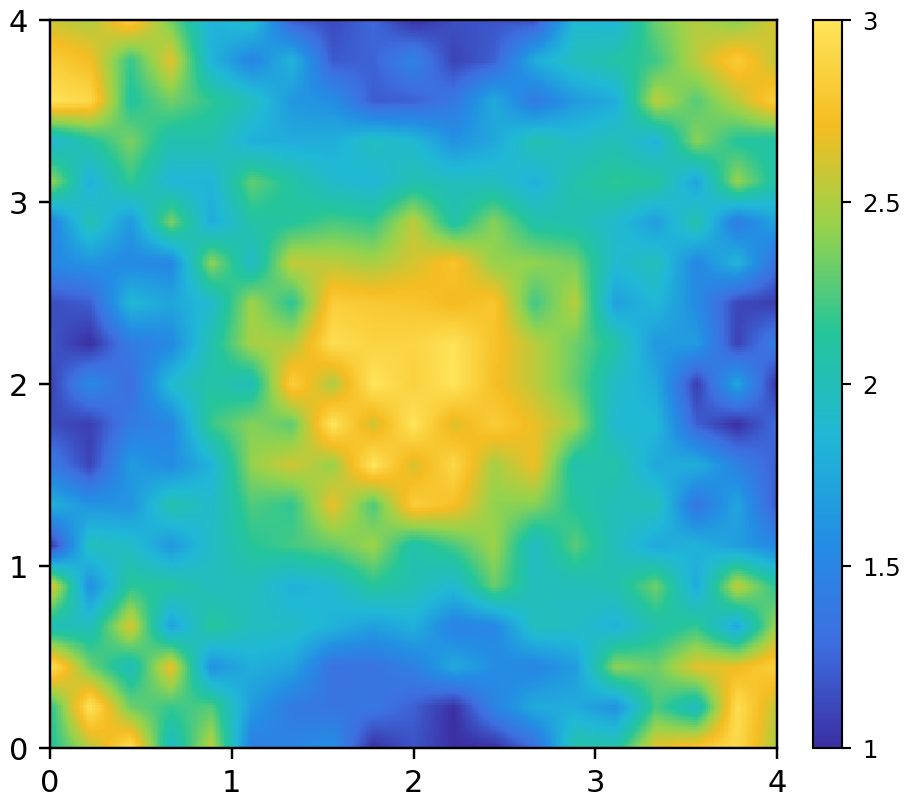}
		\end{subfigure}
		\hspace{0.025\textwidth}
		\begin{subfigure}[t]{0.23\textwidth}
			\centering
			\caption*{(b) noise 0.1\%}
			\includegraphics[width=\linewidth]{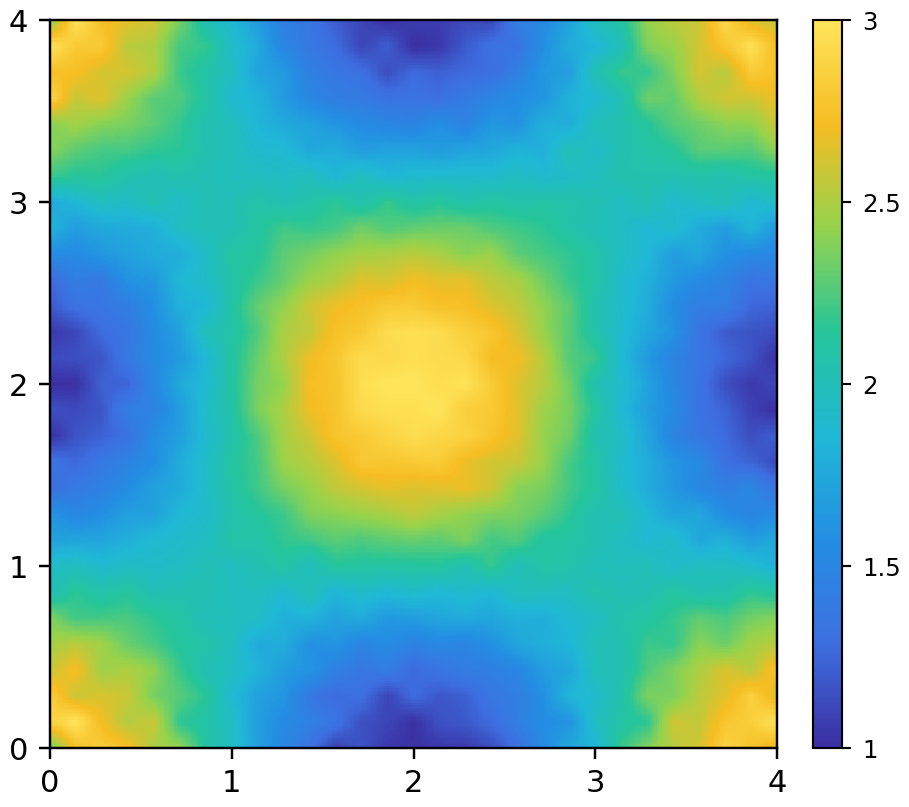}
		\end{subfigure}
		\hspace{0.025\textwidth}
		\begin{subfigure}[t]{0.23\textwidth}
			\centering
			\caption*{(c) noise 0.01\%}
			\includegraphics[width=\linewidth]{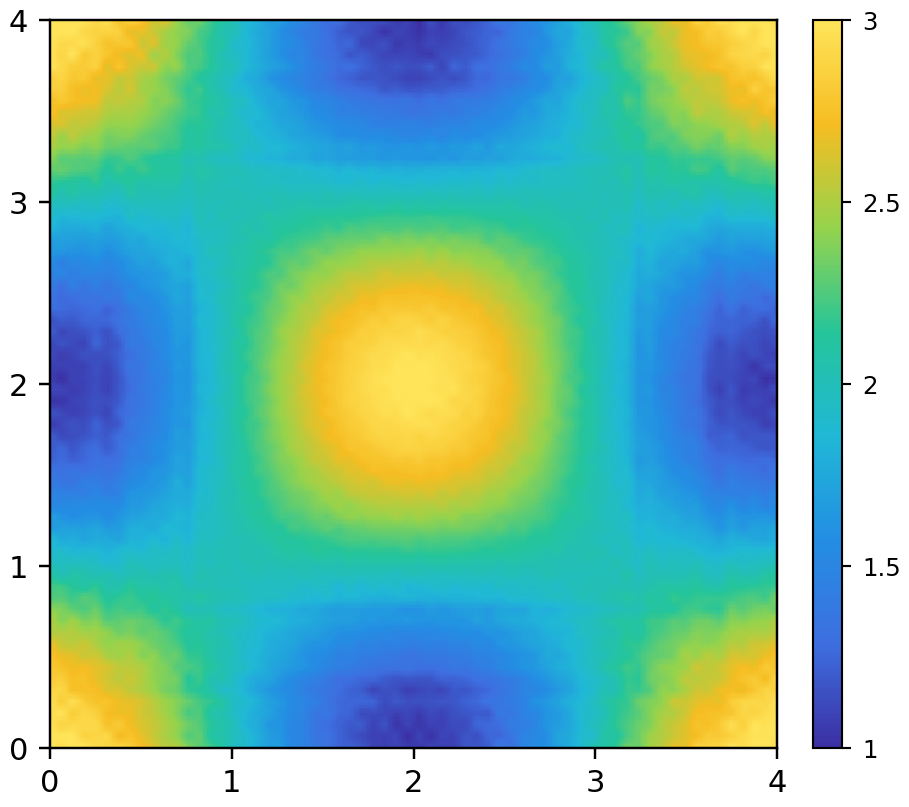}
		\end{subfigure}
		\vspace{0.25cm}
		\begin{subfigure}[t]{0.23\textwidth}
			\centering
			\includegraphics[width=\linewidth]{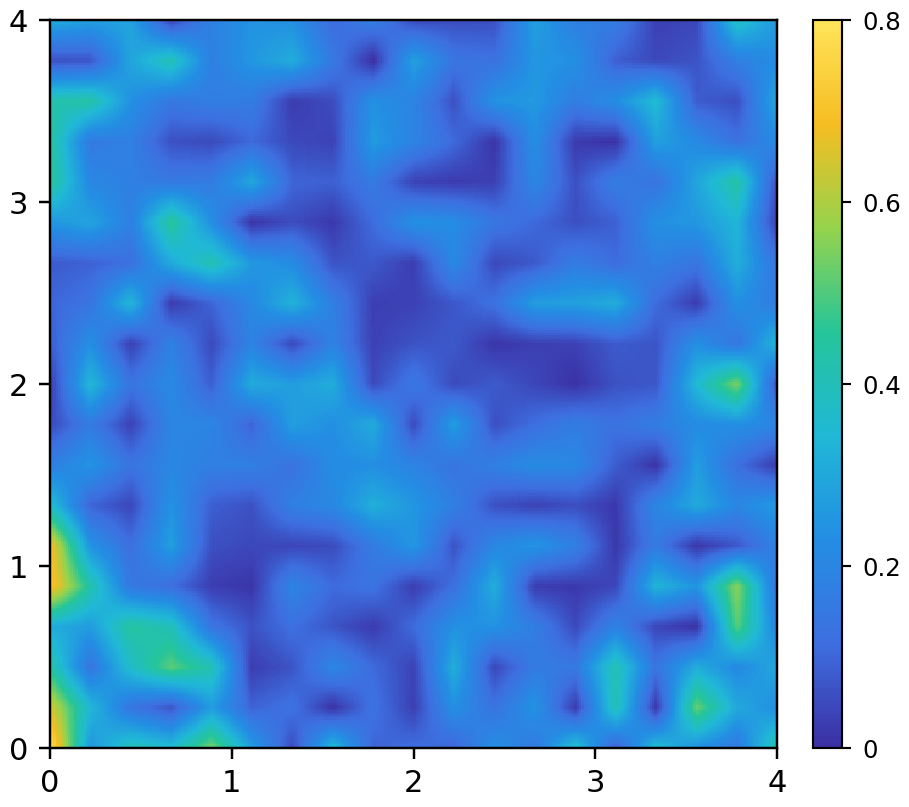}
		\end{subfigure}
		\hspace{0.025\textwidth}
		\begin{subfigure}[t]{0.23\textwidth}
			\centering
			\includegraphics[width=\linewidth]{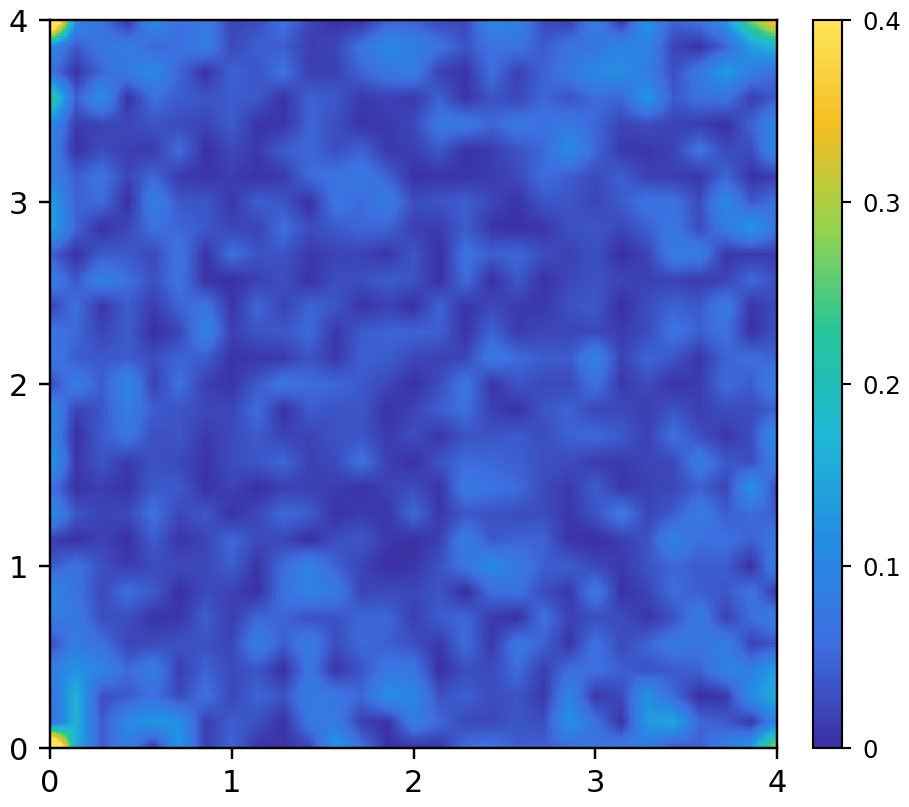}
		\end{subfigure}
		\hspace{0.025\textwidth}
		\begin{subfigure}[t]{0.23\textwidth}
			\centering
			\includegraphics[width=\linewidth]{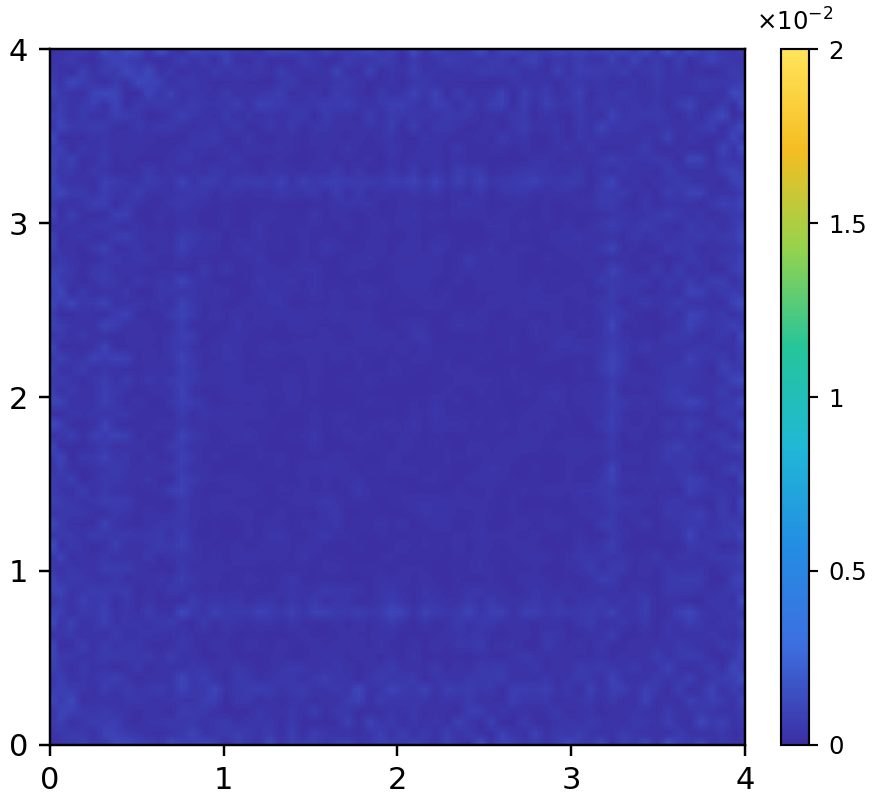}
		\end{subfigure}
		\caption{Reconstructions of smooth source \(q_1^*\) in Example~\ref{ex1}. The first row shows numerical reconstructions \(q^\delta\), and the second row shows the corresponding pointwise error \(|q^\delta-q_1^*|\).}
		\label{fig:recon-q1}
	\end{figure}

\begin{figure}
		\centering
		\begin{subfigure}[t]{0.23\textwidth}
			\centering
			\caption*{(a) noise 1\%}
			\includegraphics[width=\linewidth]{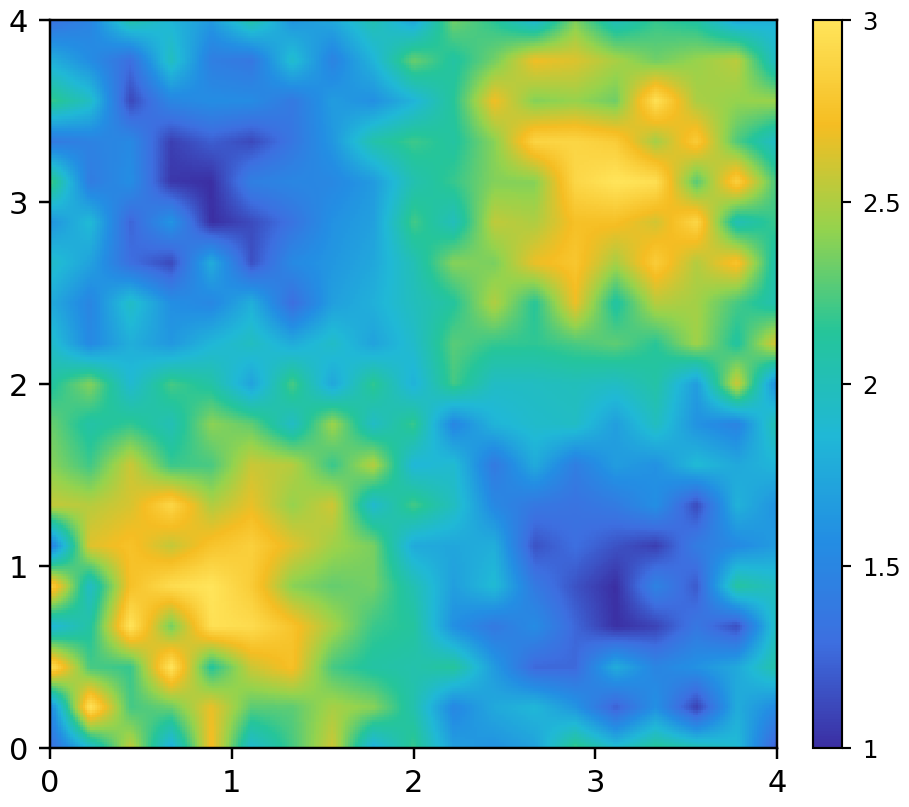}
		\end{subfigure}
		\hspace{0.025\textwidth}
		\begin{subfigure}[t]{0.23\textwidth}
			\centering
			\caption*{(b) noise 0.1\%}
			\includegraphics[width=\linewidth]{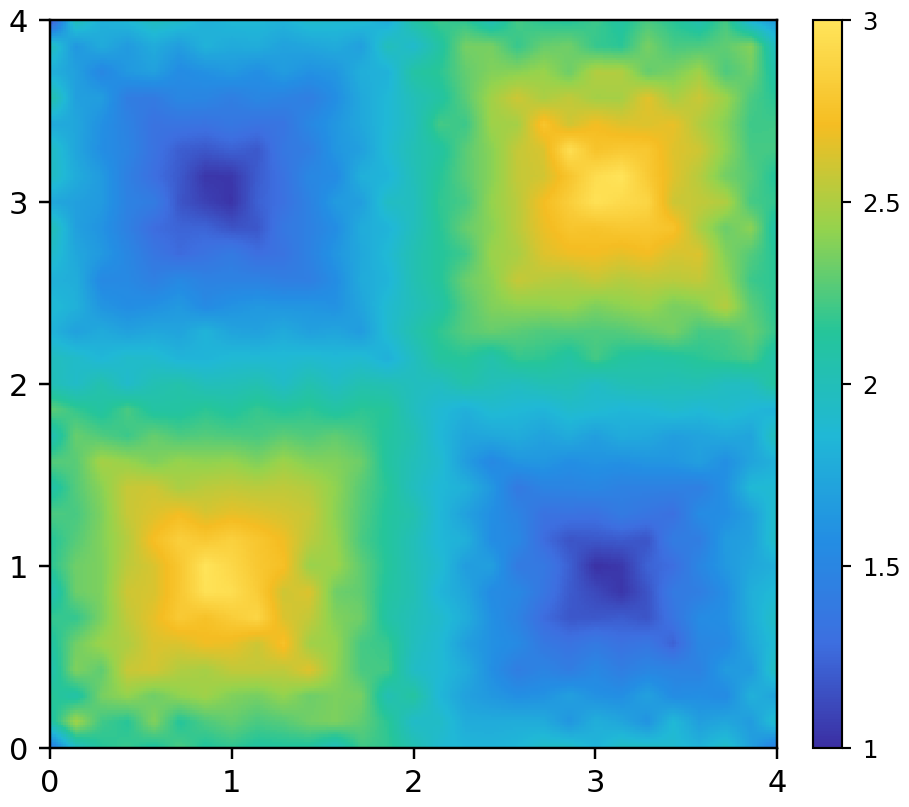}
		\end{subfigure}
		\hspace{0.025\textwidth}
		\begin{subfigure}[t]{0.23\textwidth}
			\centering
			\caption*{(c) noise 0.01\%}
			\includegraphics[width=\linewidth]{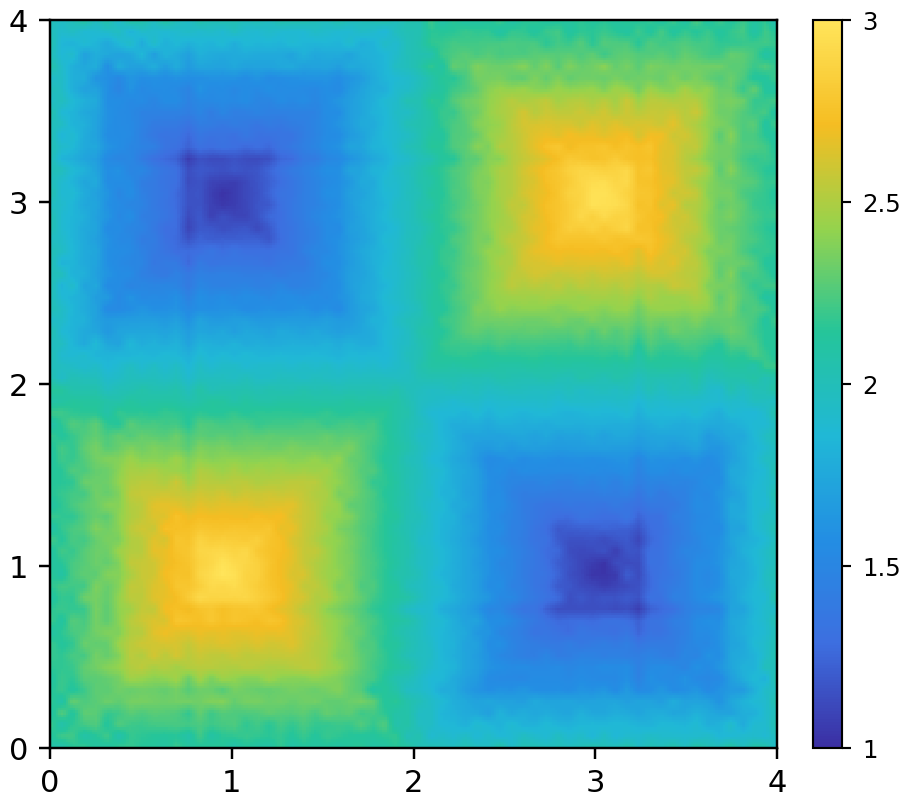}
		\end{subfigure}
		\vspace{0.25cm}
		\begin{subfigure}[t]{0.23\textwidth}
			\centering
			\includegraphics[width=\linewidth]{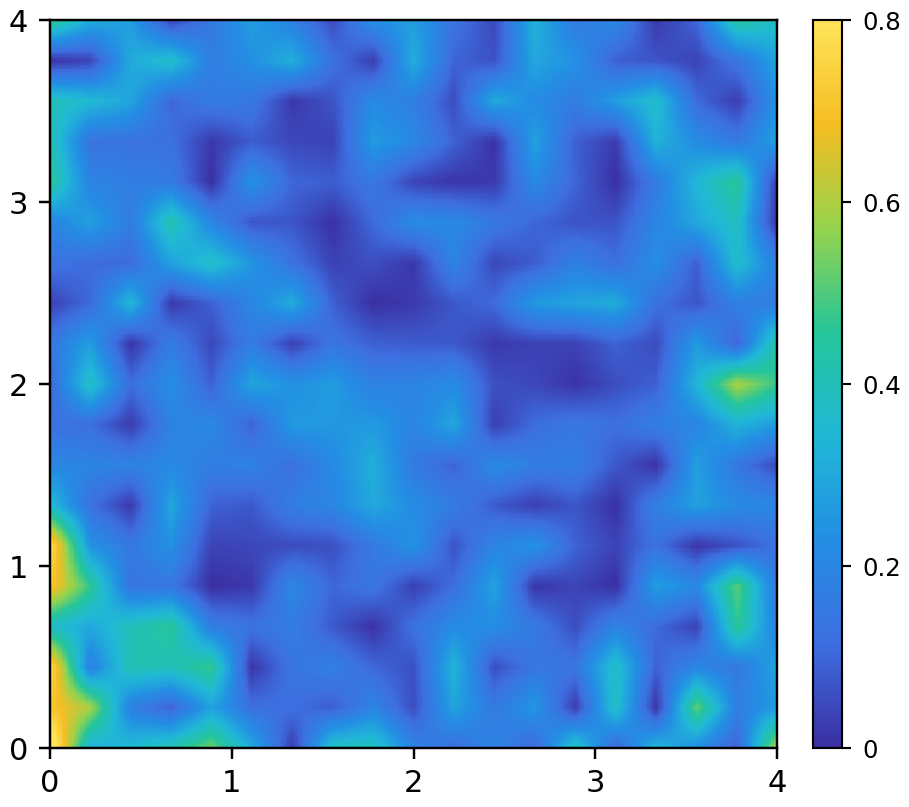}
		\end{subfigure}
		\hspace{0.025\textwidth}
		\begin{subfigure}[t]{0.23\textwidth}
			\centering
			\includegraphics[width=\linewidth]{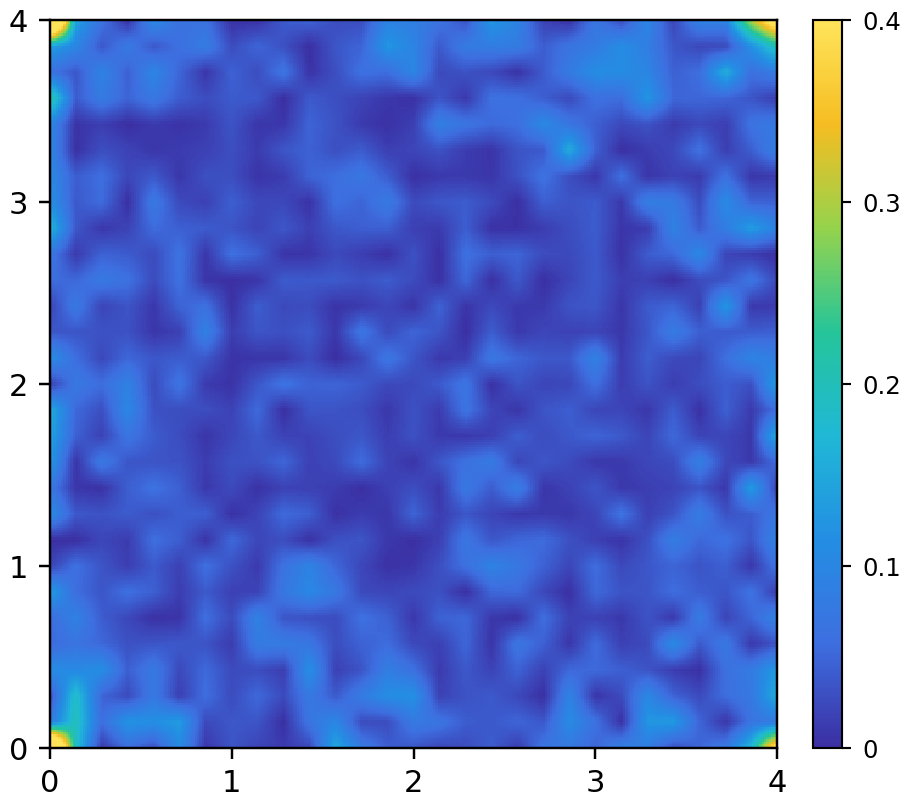}
		\end{subfigure}
		\hspace{0.025\textwidth}
		\begin{subfigure}[t]{0.23\textwidth}
			\centering
			\includegraphics[width=\linewidth]{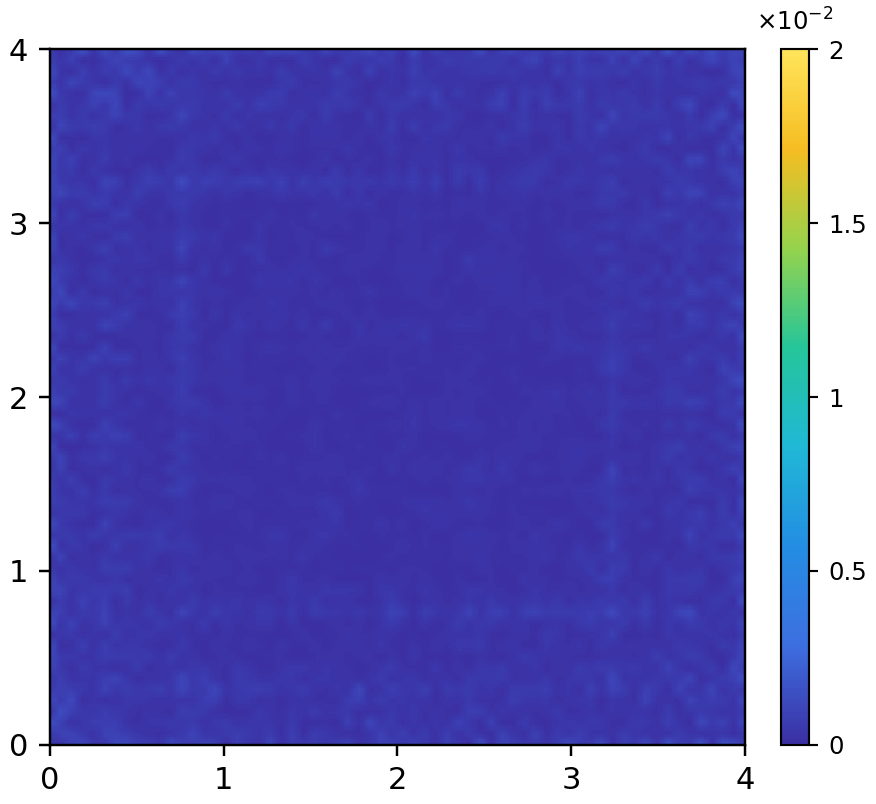}
		\end{subfigure}
		\caption{Reconstructions of piecewise smooth source \(q_2^*\) in Example~\ref{ex1}. The first row shows numerical reconstructions \(q^\delta\), and the second row shows the corresponding pointwise error \(|q^\delta-q_2^*|\).}
		\label{fig:recon-q2}
	\end{figure}

\begin{figure}
		\centering
		\begin{subfigure}[t]{0.23\textwidth}
			\centering
			\caption*{(a) noise 1\%}
			\includegraphics[width=\linewidth]{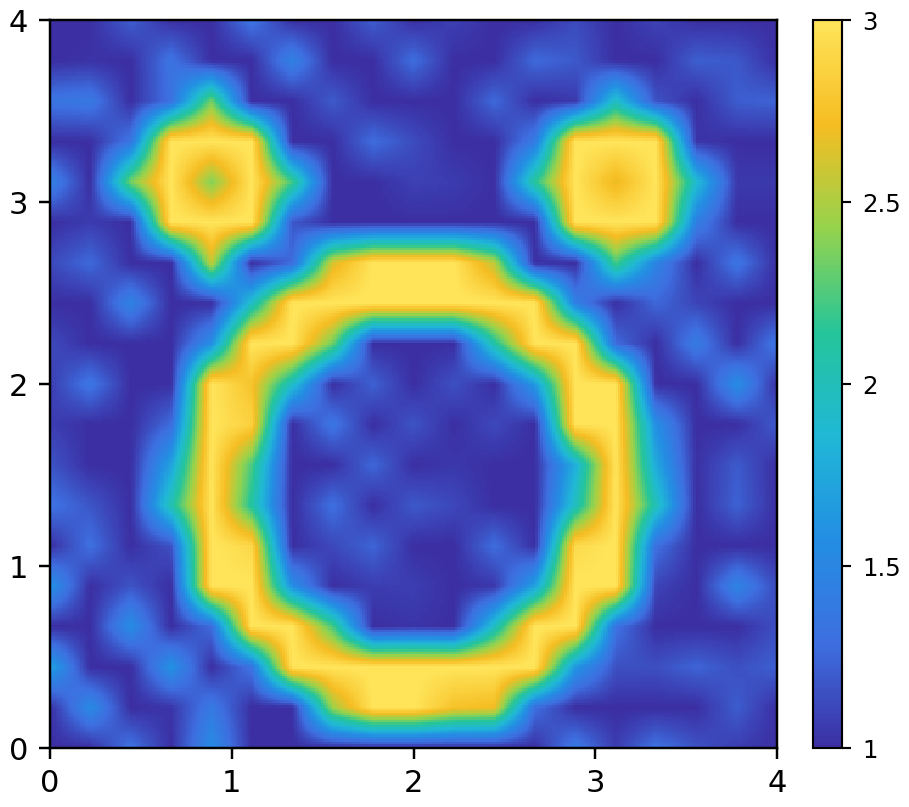}
		\end{subfigure}
		\hspace{0.025\textwidth}
		\begin{subfigure}[t]{0.23\textwidth}
			\centering
			\caption*{(b) noise 0.1\%}
			\includegraphics[width=\linewidth]{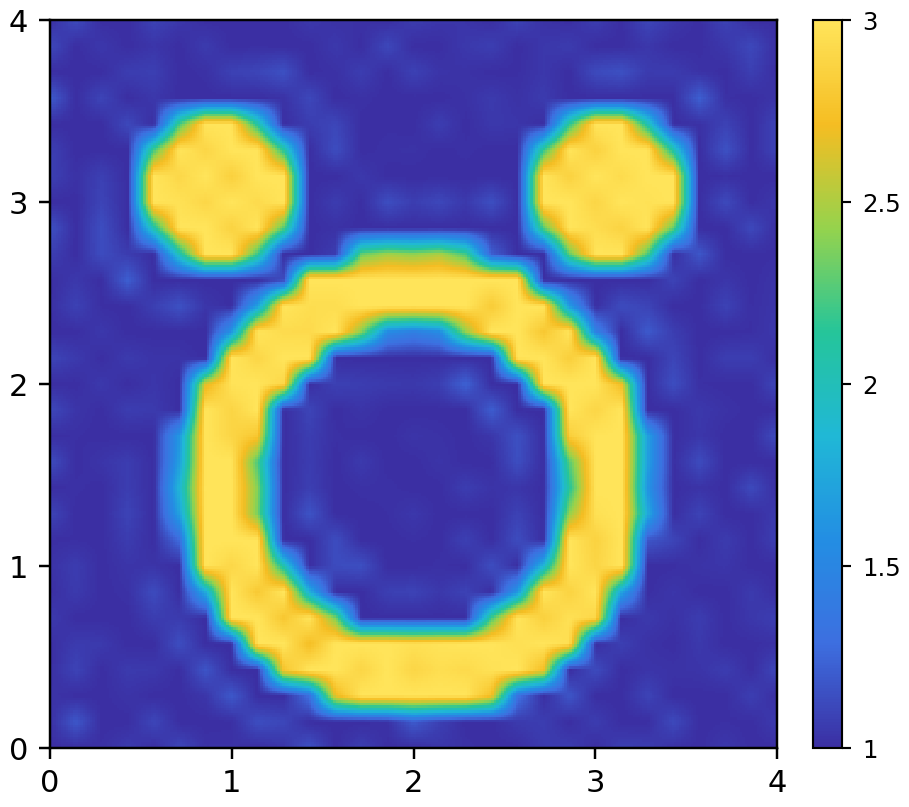}
		\end{subfigure}
		\hspace{0.025\textwidth}
		\begin{subfigure}[t]{0.23\textwidth}
			\centering
			\caption*{(c) noise 0.01\%}
			\includegraphics[width=\linewidth]{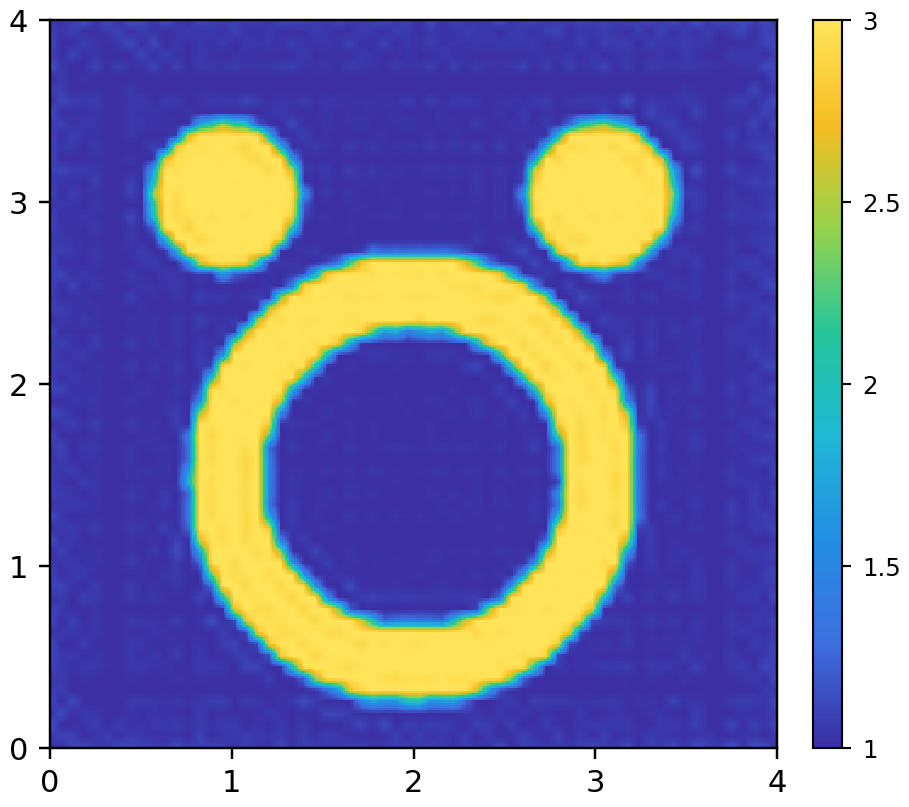}
		\end{subfigure}
		\vspace{0.25cm}
		\begin{subfigure}[t]{0.23\textwidth}
			\centering
			\includegraphics[width=\linewidth]{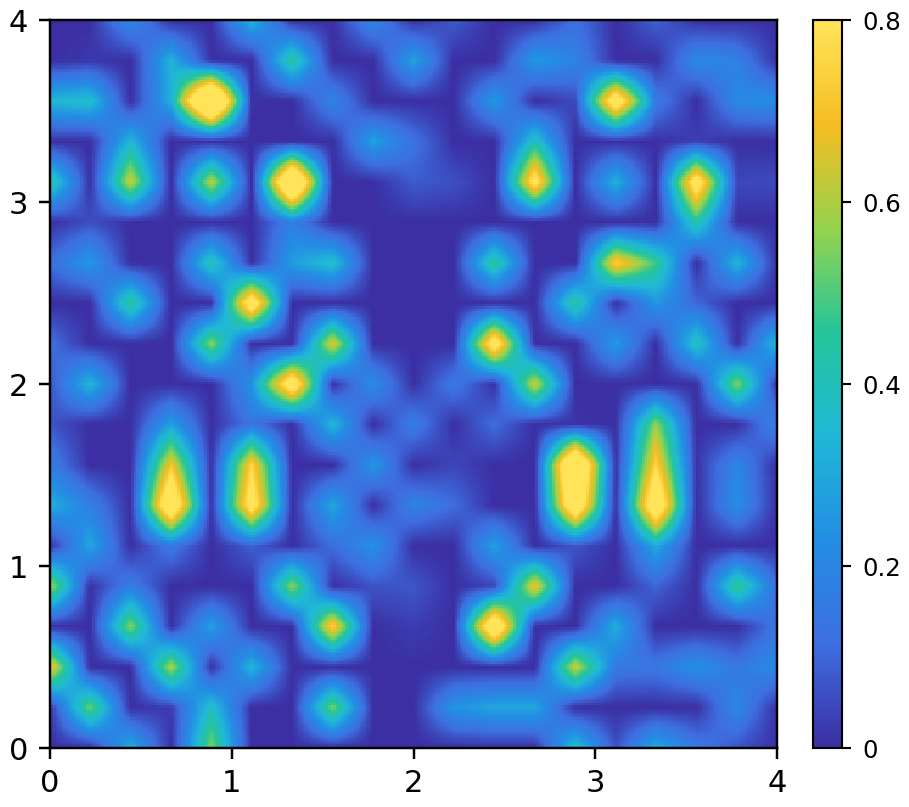}
		\end{subfigure}
		\hspace{0.025\textwidth}
		\begin{subfigure}[t]{0.23\textwidth}
			\centering
			\includegraphics[width=\linewidth]{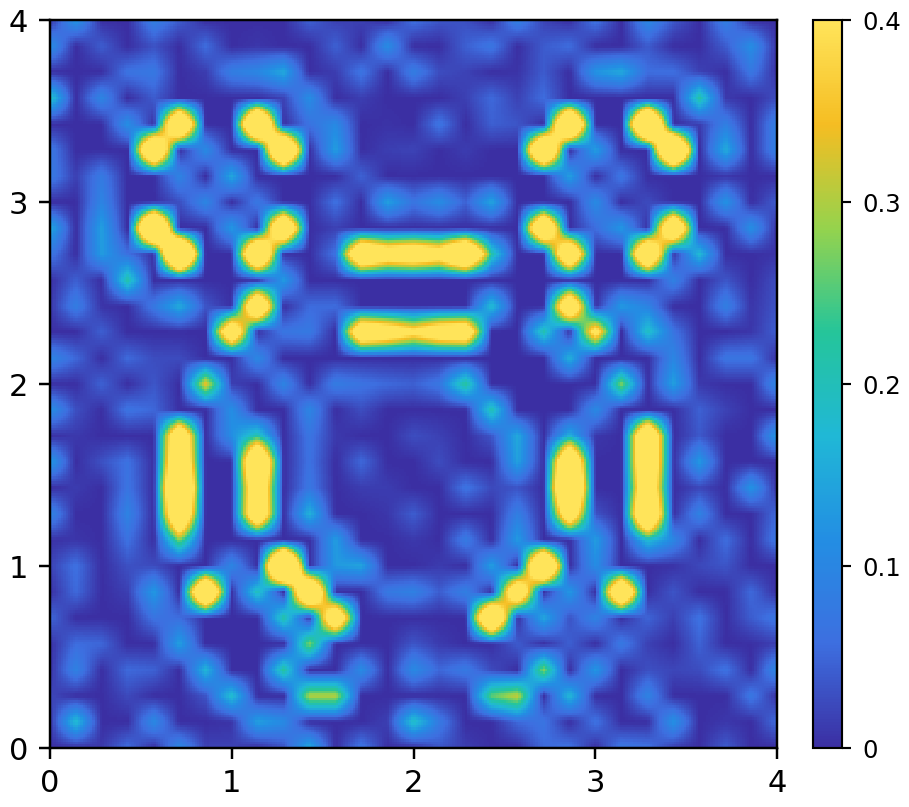}
		\end{subfigure}
		\hspace{0.025\textwidth}
		\begin{subfigure}[t]{0.23\textwidth}
			\centering
			\includegraphics[width=\linewidth]{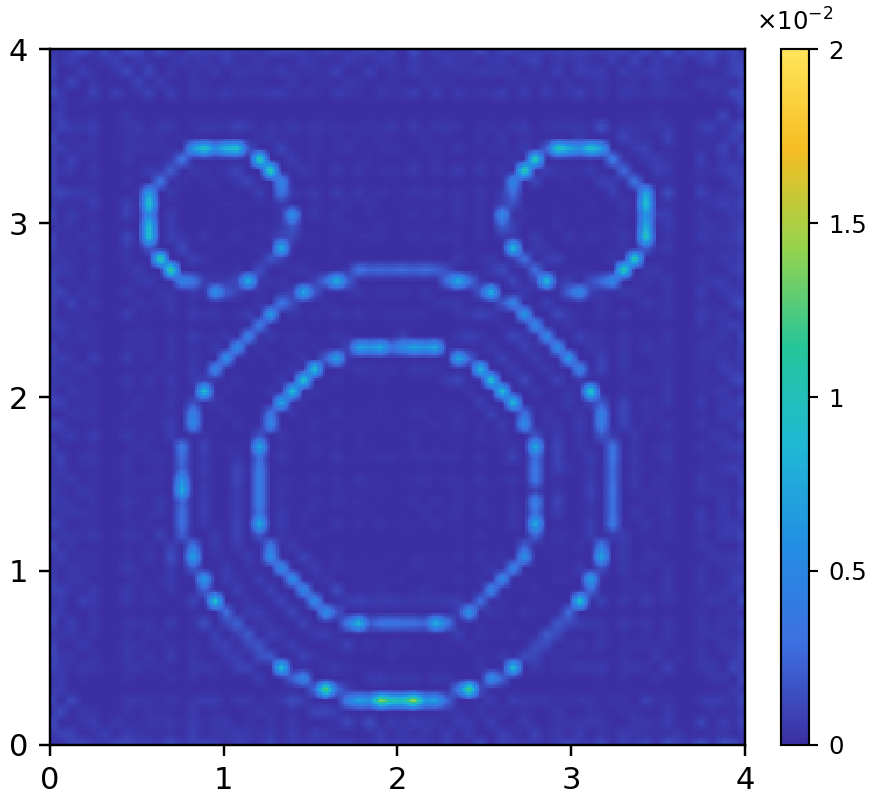}
		\end{subfigure}
		\caption{Reconstructions of discontinuous source \(q_3^*\) in Example~\ref{ex1}. The first row shows numerical reconstructions \(q^\delta\), and the second row shows the corresponding pointwise error \(|q^\delta-q_3^*|\).}
		\label{fig:recon-q3}
	\end{figure}

	Figures~\ref{fig:recon-q1}--\ref{fig:recon-q3} display the reconstructed potentials and the corresponding pointwise absolute errors for the three test potentials under different noise levels. The first row of each figure shows the numerical reconstruction \(q^\delta\), while the second row shows the pointwise error \(|q^\delta-q^*|\). As the noise level decreases, the oscillatory patterns in the reconstruction are gradually suppressed, and the pointwise error becomes smaller. For the smooth potential \(q_1^*\), the reconstruction is already quite accurate at \(1\%\) noise and becomes very close to the exact potential at \(0.01\%\) noise. For the piecewise smooth potential \(q_2^*\), the main reconstruction error is concentrated near the interfaces where the derivatives of \(q_2^*\) change. For the discontinuous potential \(q_3^*\), pronounced oscillations persist near the jump interfaces even at \(0.01\%\) noise, which is a typical feature of inverse problems with discontinuous coefficients.
	
	%%%%%%%%%%%%%%%%%%%% Figure 2 %%%%%%%%%%%%%%%%%%%
	
	%%%%%%%%%%%%%%%%%%%%%%%%% Figure 3 %%%%%%%%%%%%%%%%%
	
	%%%%%%%%%%%%%%%%%%%%%%%%% Figure 4 %%%%%%%%%%%%%%%%%
	
\begin{figure}
		\centering
		\begin{subfigure}[t]{0.305\textwidth}
			\centering
			\includegraphics[width=\linewidth]{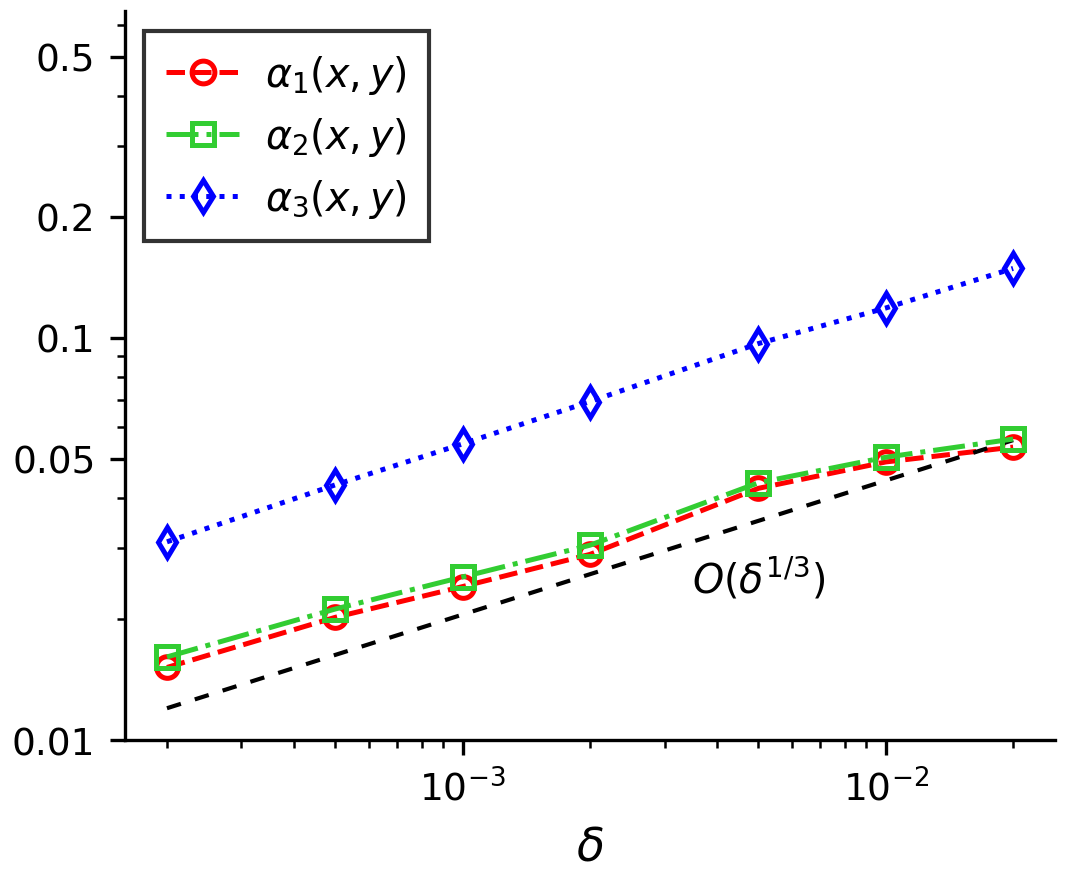}
			\caption{\(q_1^*\)}
		\end{subfigure}
		\hspace{0.01\textwidth}
		\begin{subfigure}[t]{0.305\textwidth}
			\centering
			\includegraphics[width=\linewidth]{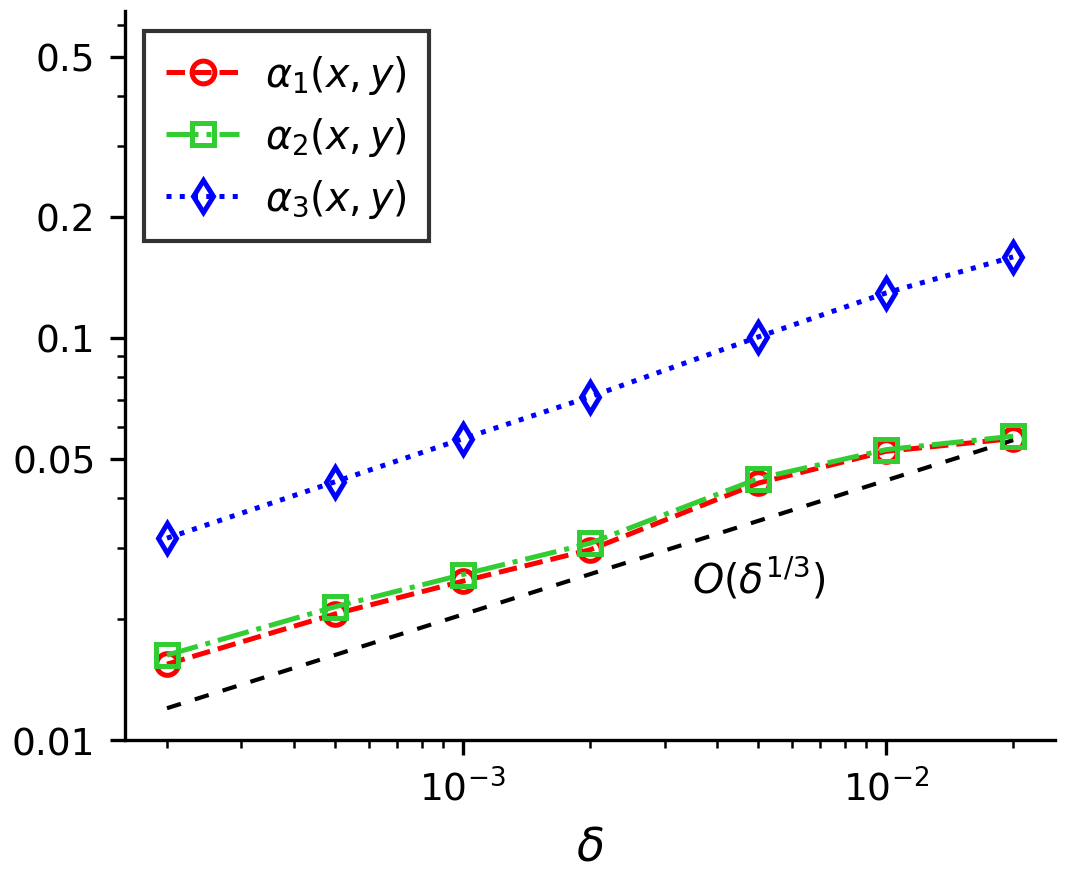}
			\caption{\(q_2^*\)}
		\end{subfigure}
		\hspace{0.01\textwidth}
		\begin{subfigure}[t]{0.315\textwidth}
			\centering
			\includegraphics[width=\linewidth]{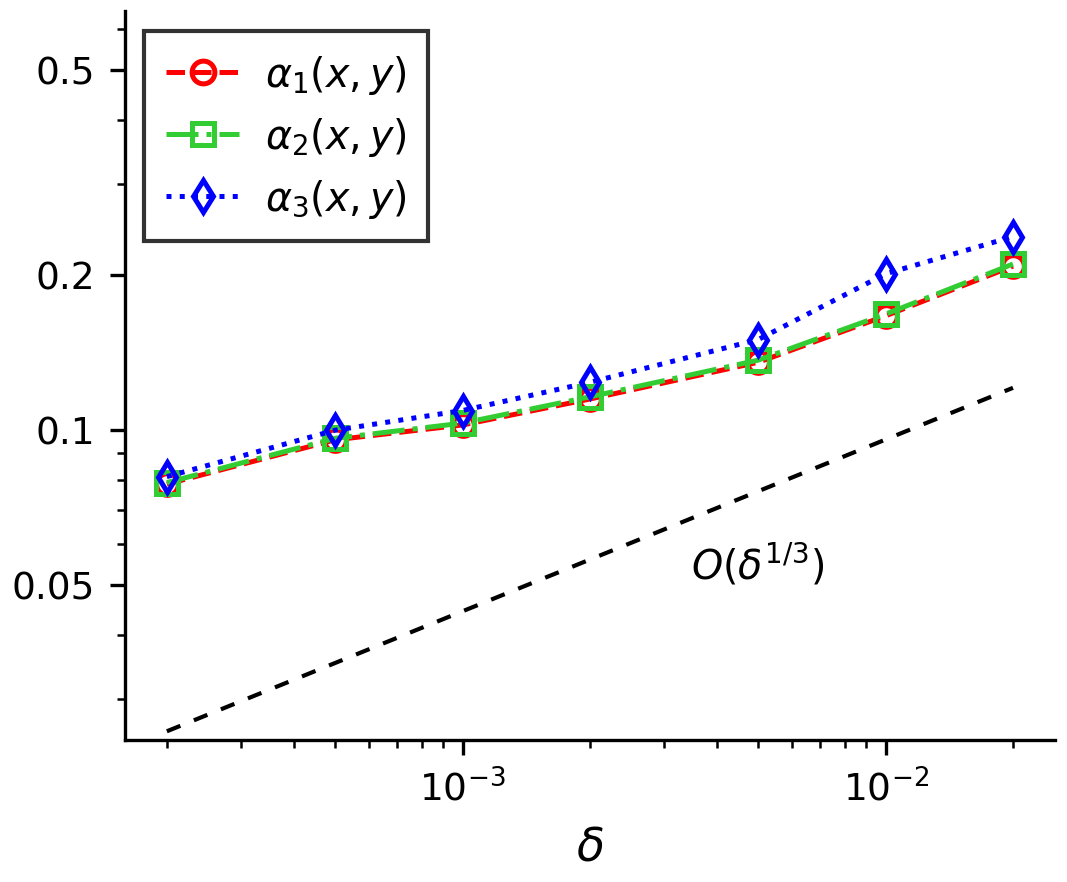}
			\caption{\(q_3^*\)}
		\end{subfigure}
		\caption{Relative error \(Err\) versus noise level \(\delta\), where \(T=5\), \(h=\frac{4}{3}\delta^{1/3}\), and \(\tau=h/15\).}
		\label{fig:convergence-delta}
	\end{figure}

	To verify the predicted rate \(\mathcal O(\delta^{1/3})\), we plot in Figure~\ref{fig:convergence-delta} the relative error \(Err\) versus the noise level \(\delta\) in a log-log scale. For the smooth potential \(q_1^*\) and the piecewise smooth potential \(q_2^*\), the observed slopes are approximately \(1/3\), which agrees well with the error estimate in Theorem~\ref{thm:inverse-error}. This confirms that the choice \(h\sim\tau\sim\delta^{1/3}\) balances the noise propagation term \(\delta/h^2\), the spatial discretization error \(h\), and the temporal discretization error \(\tau\).
	
	For the discontinuous potential \(q_3^*\), the observed convergence rate is lower than \(1/3\). This is not contradictory to the theory, since the error estimate requires sufficient regularity of \(q^*\). The numerical results therefore confirm that the regularity assumption on the exact potential is essential for achieving the optimal rate.
	
	%%%%%%%%%%%%%%%%%%%% Figure 5 %%%%%%%%%%%%%%%%%%%%%% 
\begin{figure}
		\centering
		\begin{subfigure}[t]{0.305\textwidth}
			\centering
			\includegraphics[width=\linewidth]{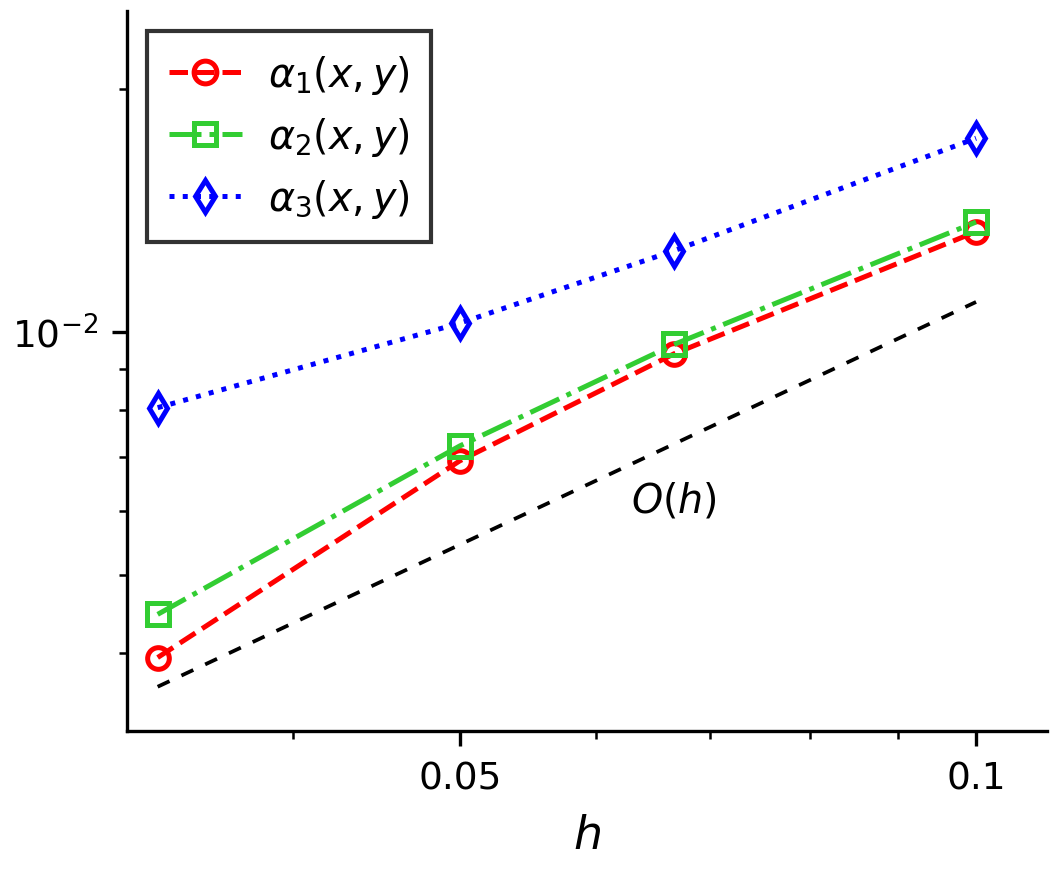}
			\caption{\(q_1^*\)}
		\end{subfigure}
		\hspace{0.01\textwidth}
		\begin{subfigure}[t]{0.305\textwidth}
			\centering
			\includegraphics[width=\linewidth]{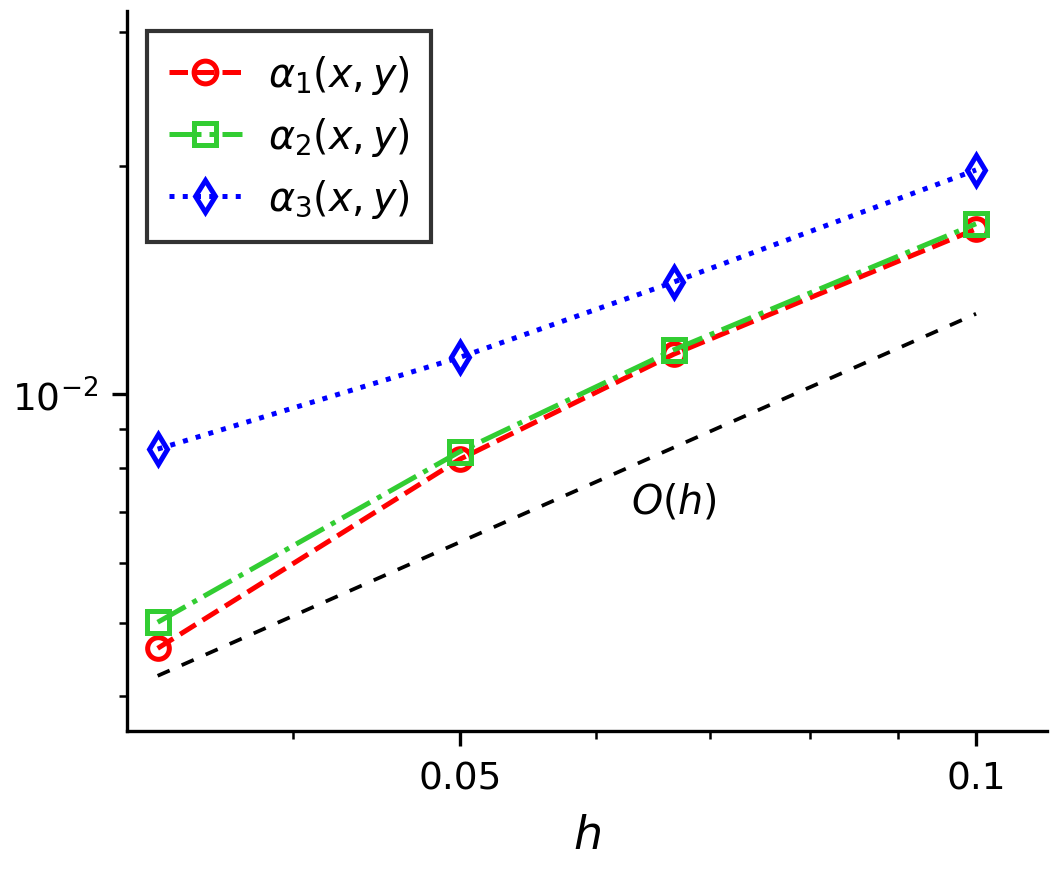}
			\caption{\(q_2^*\)}
		\end{subfigure}
		\hspace{0.01\textwidth}
		\begin{subfigure}[t]{0.315\textwidth}
			\centering
			\includegraphics[width=\linewidth]{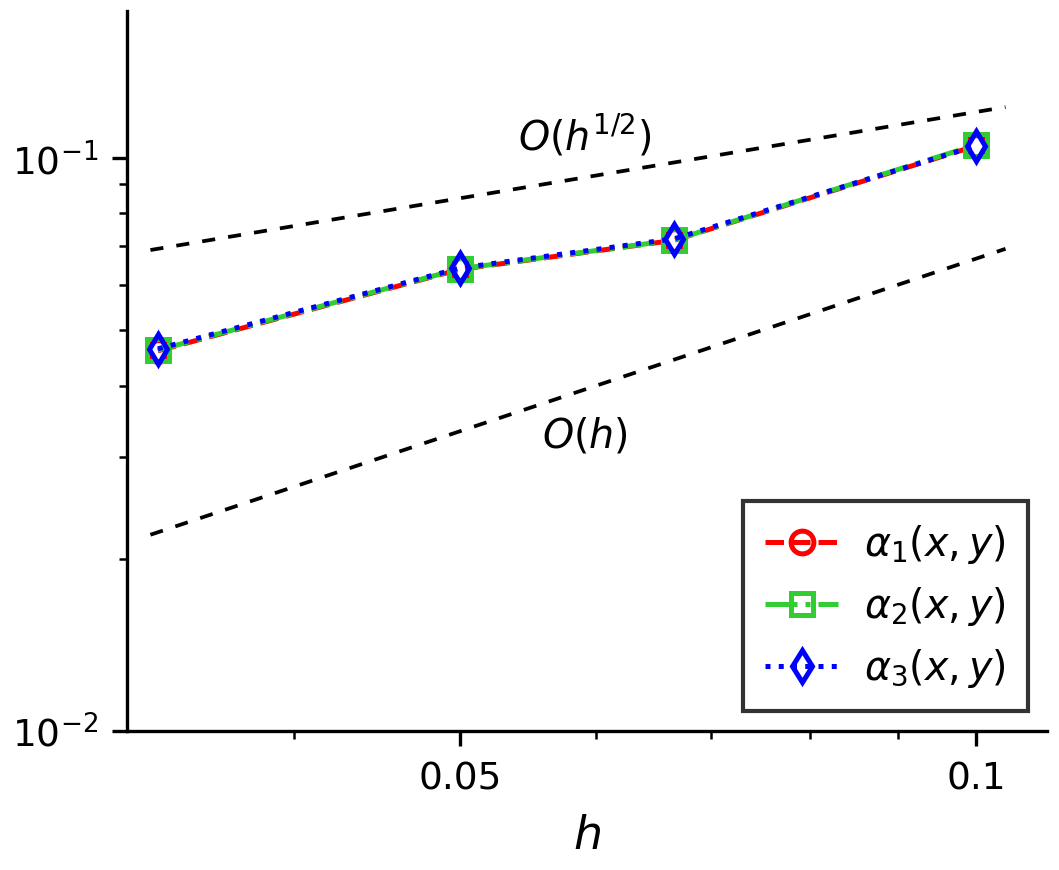}
			\caption{\(q_3^*\)}
		\end{subfigure}
		\caption{Relative error \(Err\) versus \(h\), where \(\delta=0\), \(T=5\), and \(\tau=1/500\).}
		\label{fig:convergence-h}
	\end{figure}

	To further validate the error estimate, we consider the noise-free case \(\delta=0\). In this case, the reconstruction error is dominated by the discretization error, and Theorem~\ref{thm:inverse-error} predicts $\|q^*-q^\delta\|_{L^2(\Omega)}
	=
	\mathcal O(h+\tau).$ In Figure~\ref{fig:convergence-h}, we fix \(\tau=1/500\) and plot the relative error versus the spatial mesh size \(h\). The numerical results show a clear convergence trend as \(h\) decreases, which is consistent with the spatial error term \(h\) in the theoretical estimate. The empirical convergence rate is of order $\mathcal O(h)$ for potentials $q_1^*$ and $q_2^*$. For $q_3^*$, the empirical convergence rate is of order around $\mathcal O(h^{\frac{1}{2}})$. This is due to the nonsmoothness of the discontinuous potential. These results indicate that the fully discrete scheme and the reconstruction algorithm are able to achieve the accuracy predicted by the theory when the discretization parameters are chosen appropriately.
	
	%%%%%%%%%%%%%%%%%%%%%%%%% Figure 6 %%%%%%%%%%%%%%%%%%%%
	
\begin{figure}
		\centering
		\begin{subfigure}[t]{0.305\textwidth}
			\centering
			\includegraphics[width=\linewidth]{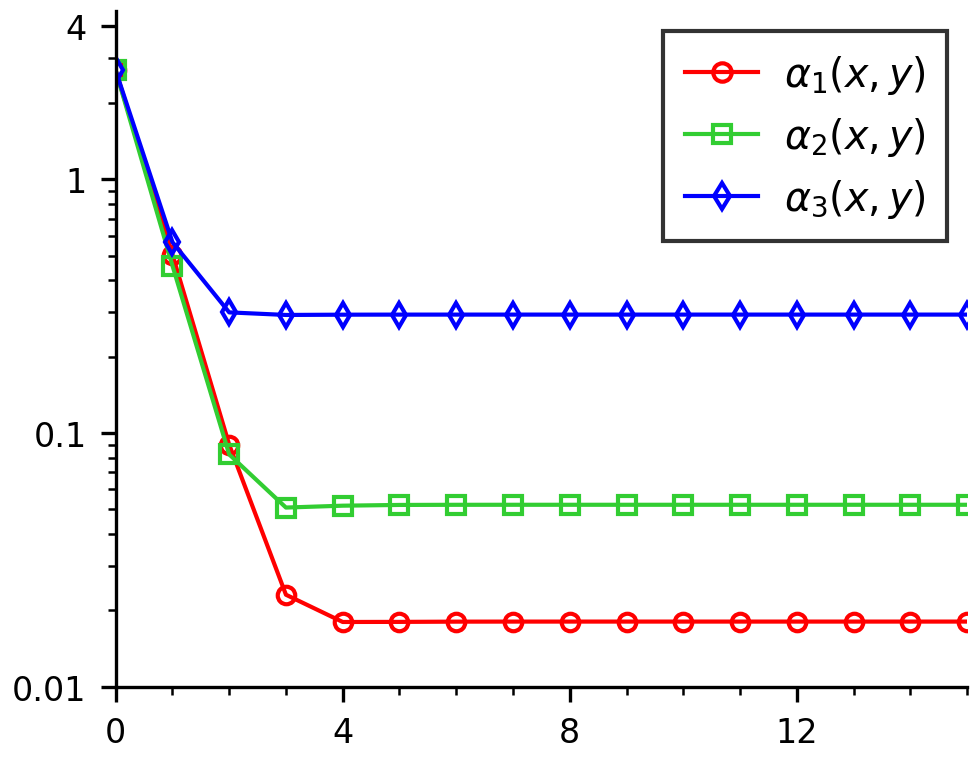}
			\caption{\(T=1\)}
		\end{subfigure}
		\hspace{0.01\textwidth}
		\begin{subfigure}[t]{0.305\textwidth}
			\centering
			\includegraphics[width=\linewidth]{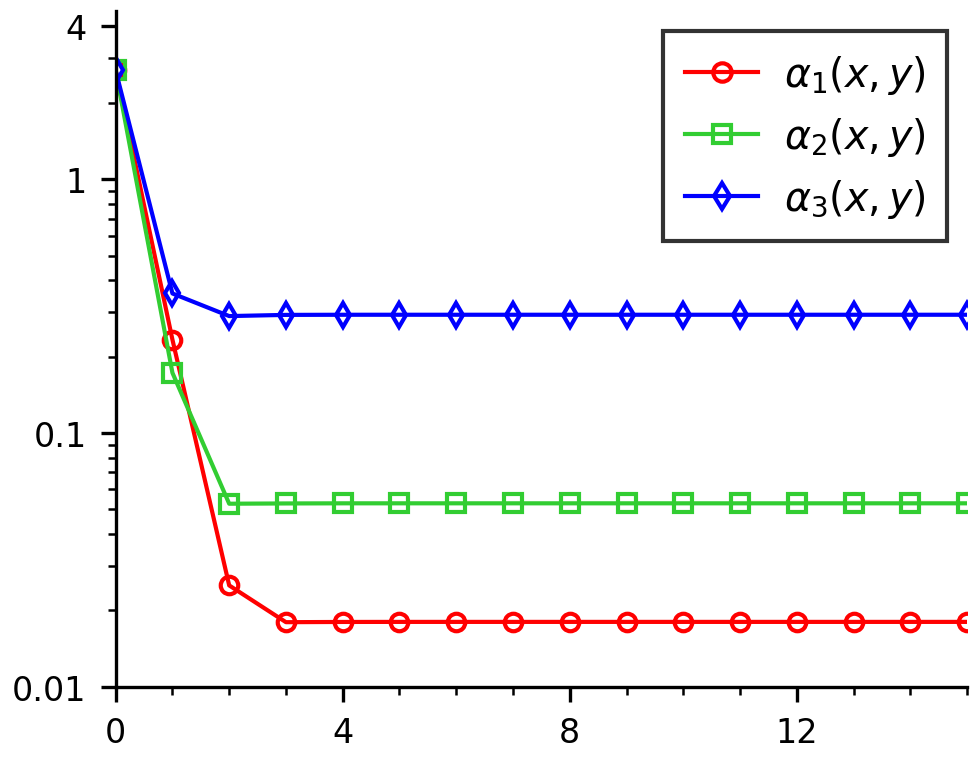}
			\caption{\(T=5\)}
		\end{subfigure}
		\hspace{0.01\textwidth}
		\begin{subfigure}[t]{0.315\textwidth}
			\centering
			\includegraphics[width=\linewidth]{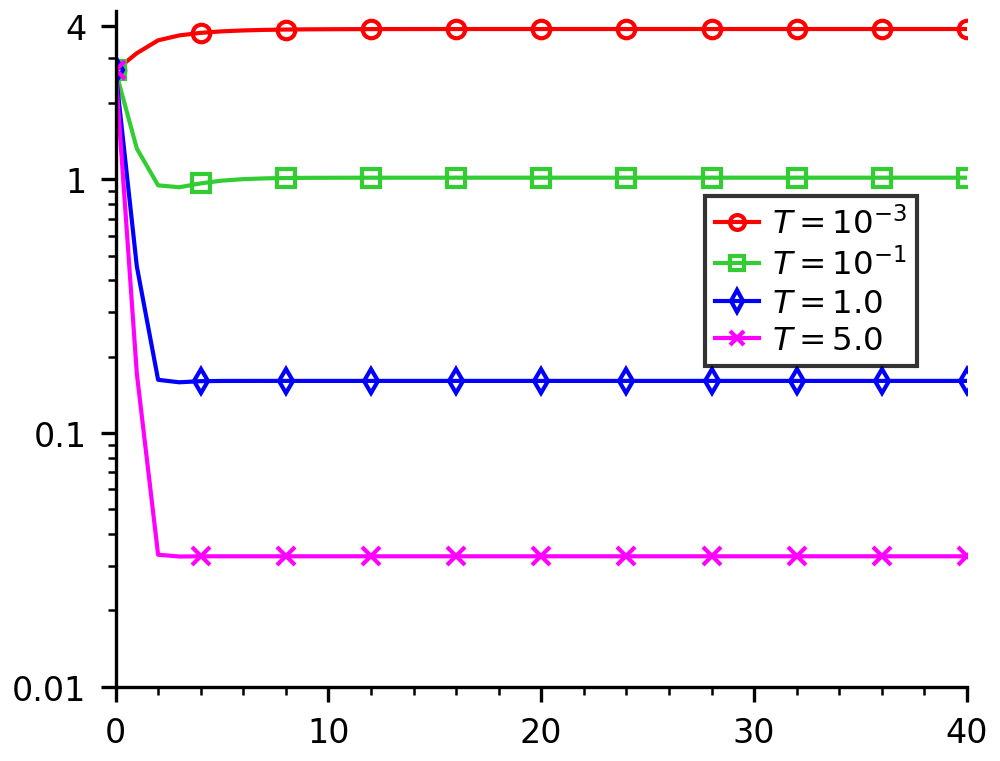}
			\caption{\(\alpha=\alpha_2\)}
		\end{subfigure}
		\caption{Convergence histories of Algorithm~\ref{alg:recon} with different \(T\) and \(\alpha\), where (a)-(b): \(\delta=10^{-6}\), \(h=0.1\), and \(\tau=h/5\), (c): \(\delta=10^{-6}\), \(h=0.1\), and \(\tau=T/150\).}
		\label{fig:convergence-T-alpha}
	\end{figure}

	We investigate the influence of the terminal time \(T\) and the variable exponent \(\alpha(x)\) on the convergence of Algorithm~\ref{alg:recon}. In Figure~\ref{fig:convergence-T-alpha}, we report the convergence histories of the iteration for different \(T\) and \(\alpha\). The vertical axis is the error between the current iteration and the exact potential, and the horizontal axis is the iteration index. The curves exhibit a linear decay in a semi-logarithmic scale, which confirms the linear convergence of the fixed-point iteration established in Lemma~\ref{lem:contraction-Kdelta} and Theorem~\ref{thm:inverse-error}.
	
	Moreover, the convergence factor decreases as \(T\) increases. This agrees with the contraction estimate $\|K_{h,\tau}^{\delta}q_1-K_{h,\tau}^{\delta}q_2\|_{L^2(\Omega)}
	\le
	C\max\{T^{-\alpha_*},T^{-(1-\epsilon)\alpha_*}\}
	\|q_1-q_2\|_{L^2(\Omega)},$ which shows that a larger terminal time leads to a smaller contraction factor and hence faster convergence. The experiments also indicate that the iteration is robust with respect to the variable exponent \(\alpha(x)\): although the value of \(\alpha\) affects the constants and the decay rate of the fractional derivative, the observed convergence behavior remains stable. This is consistent with the theoretical analysis, where the variable exponent enters the estimates only through \(\alpha_*\) and \(\alpha^*\).
	
	%%%%%%%%%%%%%%%%%%%%%%% Figure 7 %%%%%%%%%%%%%%%%%%%%%

	The condition \(T\ge T_*\) in Theorems~\ref{thm-stability} and~\ref{thm:inverse-error} is not merely technical. To illustrate this, we also tested the algorithm with very small terminal times, e.g., \(T=10^{-1}\) and \(T=10^{-3}\) in Figure \ref{fig:convergence-T-alpha}(c). In these cases, the terminal observation contains insufficient information about the potential, and the contraction factor in Lemma~\ref{lem:contraction-Kdelta} is close to or larger than one. As a result, the iteration either stagnates or converges to a limit far from the exact potential. This confirms that a sufficiently large terminal time is necessary for stable and accurate reconstruction, as predicted by the theory.

	%\paragraph{Summary.}
	%In summary, the numerical experiments support the theoretical results in the following aspects:
	%\begin{itemize}
	%\item the fixed-point iteration converges linearly, as predicted by the contraction estimate;
	%\item for smooth and piecewise smooth potentials, the reconstruction error exhibits the predicted rate \(\mathcal O(\delta^{1/3})\) when \(h\sim\tau\sim\delta^{1/3}\);
	%\item the spatial and temporal discretization errors behave like \(\mathcal O(h+\tau)\), which validates the error estimate in Theorem~\ref{thm:inverse-error};
	%\item the reconstruction is stable and accurate for sufficiently large terminal time \(T\);
	%\item the regularity of the exact potential is important: discontinuous potentials lead to larger errors and slower convergence rates;
	%\item the variable order \(\alpha(x)\) affects the constants but not the overall convergence mechanism, which confirms the robustness of the proposed method for variable-exponent sub-diffusion.
	%\end{itemize}

	\section{Conclusion}
	
	In this paper, we have investigated an inverse potential problem for a variable-exponent sub-diffusion equation from terminal observation data. We formulated the inverse problem as a nonlinear fixed-point equation and proved that the associated operator \(K\) is monotone, which leads to a decreasing fixed-point iteration converging to the unique solution. The uniqueness and conditional stability of the inverse problem were established under suitable assumptions on the problem data and a sufficiently large terminal time \(T\). For the numerical reconstruction, we developed a fully discrete scheme based on the finite element method in space and a convolution quadrature in time, and we proved the linear convergence of the iterative algorithm. An a priori error estimate was derived, which provides a guideline for choosing the discretization parameters according to the noise level. Specifically, the choice \(h\sim\tau\sim\delta^{1/3}\) yields the optimal convergence rate \(\mathcal O(\delta^{1/3})\).
	
	Numerical experiments in a two-dimensional setting confirmed the theoretical findings. The results showed that the fixed-point iteration converges linearly and that the reconstruction error exhibits the predicted rate for smooth and piecewise smooth potentials. The spatial and temporal discretization errors were also validated. The experiments illustrated the necessity of a sufficiently large terminal time and the regularity assumption on the exact potential. Moreover, the method was shown to be robust with respect to the variable exponent \(\alpha(x)\).
	
	\section*{Acknowledgements}
	Zhiyuan Li is supported by the National Natural Science Foundation of China (Grant No. 12271277), the Ningbo Youth Leading Talent Project (Grant No. 2024QL045), and the Open Research Fund of Key Laboratory of Nonlinear Analysis \& Applications (Central China Normal University), Ministry of Education, China. Chunlong Sun is supported by the Fundamental Research Funds for the Central Universities (Grant No. NS2026033) and the Key Laboratory of Nonlinear Analysis $\&$ Applications (Ministry of Education), Central China Normal University, Wuhan 430079, China. Xiangcheng Zheng is supported by the National Natural Science Foundation of China (Nos. 12301555 and 12622126) and the Natural Science Foundation of Shandong Province (Grant No. ZR2025QB01).

\end{document}